\documentclass[twocolumn,amsthm]{autart}

\usepackage{amsmath,amssymb,mathtools}
\usepackage{enumitem}
\usepackage{array}
\usepackage{graphicx}
\usepackage{xcolor}
\usepackage{booktabs}
\usepackage{subcaption}
\usepackage{natbib}
\usepackage[vlined,ruled,norelsize,algo2e]{algorithm2e}

\bibpunct{(}{)}{;}{a}{,}{,}

\newcommand{\calF}{\mathcal{F}}

\newcommand{\calN}{\mathcal{N}}

\newcommand{\calV}{\mathcal{V}}

\newcommand{\Chat}{\hat{C}}

\newcommand{\rbar}{\bar{r}}

\newcommand{\xbar}{\bar{x}}

\newcommand{\gammabar}{\bar{\gamma}}

\newcommand{\mrm}{\mathrm{m}}

\newcommand{\srm}{\mathrm{s}}

\newcommand{\Nrm}{\mathrm{N}}

\renewcommand{\Re}{\mathbb{R}}

\newcommand{\Ne}{\mathbb{N}}

\newcommand{\dxd}{{d\times d}}

\newcommand{\diff}{\mathrm{d}}

\newcommand{\maxtxt}{\mathrm{max}}

\newcommand{\Sbb}{\mathbb{S}}
\renewcommand{\Bbb}{\mathbb{B}}

\newcommand{\Vmax}{\calV_{\maxtxt}}
\newcommand{\amax}{a_{\maxtxt}}

\begin{document}

\begin{frontmatter}
\title{Counterexample-guided computation of polyhedral Lyapunov functions for hybrid systems\thanksref{footnoteinfo}}

\thanks[footnoteinfo]{G.\ Berger is a BAEF fellow. S.\ Sankaranarayanan was supported by the 
US National Science Fundation (NSF). }

\author[CUBoulder]{Guillaume O.\ Berger}\ead{guillaume.berger@colorado.edu},            % Add the
\author[CUBoulder]{Sriram Sankaranarayanan}\ead{sriram.sankaranarayanan@colorado.edu},  % e-mail address
                                                                                        % (ead) as shown

\address[CUBoulder]{Department of Computer Science, University of Colorado, Boulder, USA}   % Please supply
                                                                                            % full addresses
                                                                                            % here.

% Five to ten keywords, chosen from the IFAC keyword list or with the
% help of the Automatica keyword wizard

\begin{abstract} % Abstract of not more than 200 words.
This paper presents a counterexample-guided iterative algorithm to compute convex, piecewise linear (polyhedral) Lyapunov functions for uncertain continuous-time linear hybrid systems. Polyhedral Lyapunov functions provide an alternative to commonly used polynomial Lyapunov functions. Our approach first characterizes intrinsic properties of a polyhedral Lyapunov function including its ``eccentricity'' and ``robustness'' to perturbations. We then derive an algorithm that either computes a polyhedral Lyapunov function proving that the system is stable, or concludes that no polyhedral Lyapunov function exists whose eccentricity and robustness parameters satisfy some user-provided limits. Significantly, our approach places no a priori bounds on the number of linear pieces that make up the desired polyhedral Lyapunov function. 

The algorithm alternates between a learning step and a verification step, always maintaining a finite set of witness states. The learning step solves a linear program to compute a candidate Lyapunov function compatible with a finite set of witness states. In the verification step, our approach verifies whether the candidate Lyapunov function is a valid Lyapunov function for the system. If verification fails, we obtain a new witness. We prove a theoretical bound on the maximum number of iterations needed by our algorithm. We demonstrate the applicability of the algorithm on numerical examples.
\end{abstract}

\end{frontmatter}

%%%%%%%%%%%%%%%%%%%%%%%%%%%%%%%%%%%%%%%%%%%%%%%%%%%%%%%%%%%%%%%%%%%%%%%%%%%%%%%%
\section{Introduction}\label{sec-introduction}

Lyapunov methods are a very useful tool for dynamical systems 
analysis~\citep{khalil2002nonlinear,liberzon2003switching}.
The existence of a Lyapunov function for a given system allows to 
study important properties of the system, such as stability or positive invariants. 
However, the problem of finding such a function for a given system is quite 
challenging, especially for hybrid systems that switch between multiple 
modes with varying dynamics.  In this paper, 
we consider the problem of automatically synthesizing Lyapunov 
functions to prove the stability of a class of uncertain piecewise linear systems. Such 
systems can switch between various modes wherein the dynamics for each mode 
are linear.  Proving the stability 
of a piecewise linear system can be achieved using a 
\emph{common Lyapunov function} which is shown to be decreasing 
according to the dynamics of all the modes of the system~\citep{sun2011stability}. However, 
finding a common Lyapunov function can be challenging, as well. Many methods 
search for \emph{polynomial} Lyapunov functions using approaches 
such as \emph{Sum-of-Squares} (SOS) programming~\citep{lasserre2001global}. 
However, such methods  also restrict the degree of the desired Lyapunov function to be within 
some a priori bounds. Failure to find such a function gives us little insight 
as to the nature of the underlying system. 

In this paper, we focus on common Lyapunov functions that are convex and 
piecewise linear (also known as polyhedral) functions. Polyhedral functions 
are an interesting and relatively less studied class of Lyapunov functions when compared 
to polynomial Lyapunov functions.  
Their expressiveness can be  modulated by adjusting the number of linear 
pieces~\citep{blanchini2015settheoretic}. Furthermore, for a 
large class  systems, including switched linear systems, there exist inverse results that 
show that if the system is stable, then a polyhedral Lyapunov function exists~\citep{sun2011stability}.
The computation of polyhedral Lyapunov functions nevertheless remains 
difficult, even for linear systems.
The main difficulties arise from the fact that the requirement for the 
function to decrease along the trajectories of the system gives rise to 
nonlinear, nonconvex constraints~\citep{blanchini2015settheoretic}. Furthermore, 
there is no known a priori bound on the number of pieces that such a function 
must have in order to be a Lyapunov function for a given class of systems~\citep{ahmadi2016lower}. 

We introduce a  counterexample-guided approach to compute polyhedral 
Lyapunov functions for continuous-time, uncertain, piecewise linear dynamical systems.
First, we formally characterize polyhedral Lyapunov functions in terms of  key 
parameters that include  (a) the \emph{eccentricity} of the function which is defined 
analogously to the well-known concept  of eccentricity of an ellipsoid; and (b) the \emph{robustness} 
of the function which is defined in terms of limits on the perturbations of the dynamics 
of the underlying system that can be proven  stable by the given polyhedral Lyapunov function. 
Naturally, we desire Lyapunov functions  whose eccentricities are as small as possible 
and robustness parameters are as large as possible. In addition to the system being analyzed, our 
algorithm takes as inputs limits on the eccentricity and robustness parameters. It can yield 
two possible outcomes: (a) success along with a polyhedral Lyapunov function that proves the stability 
of the given system; or (b) failure to find a suitable function. However, unlike the related work
we guarantee that if our algorithm fails then \emph{no polyhedral Lyapunov function exists} for the 
system at hand, whose eccentricities and robustness parameters lie within the user-provided limits. 
By adjusting these limits, we can trade-off between the ``complexity'' of the Lyapunov functions
we are searching for against computation time and the precision needed for the numerical computations 
involved. Significantly, our algorithm does not place any \emph{a priori} bounds on the number of pieces of 
the desired Lyapunov function. 

The algorithm itself is iterative, maintaining a finite set of states that we will 
call the \emph{witness set}. Each iteration performs a learning step, followed by 
a verification step. In the learning step, a candidate polyhedral Lyapunov function is 
computed based on the witness set.  Subsequently, the verification step checks whether 
the candidate Lyapunov function is a  valid Lyapunov function for the system, and if not, 
outputs a new state wherein the Lyapunov conditions  are violated. This state is added 
to the witness set in order to ensure that a previous candidate 
is never re-visited by our algorithm. The learning--verification procedure 
is repeated, until no counterexample is found in the verification step in which case 
our current candidate is the desired Lyapunov function or the learning step fails to 
identify a candidate which allows us to conclude that no Lyapunov function exists within
user-provided eccentricity and robustness parameters. 

A desirable property of our approach is that the learning and verification steps are 
implemented by solving a series of linear programming problems whose sizes are bounded 
by those of the  underlying system and the number of witnesses found thus far. 
At the same time, we prove that the number of iterations of our approach is bounded and 
derive upper bounds in two different ways by establishing key properties of polyhedral functions 
whose eccentricities and robustness parameters are bounded. 

The use of finitely many witness states allows our algorithm to avoid bilinear (nonconvex) 
constraints that often arise in \emph{direct approaches} that try to enforce the Lyapunov conditions 
for an  unknown polyhedral function over all (infinitely many) states.  

We evaluate our approach on a series of numerical examples ranging from challenging 
instances that have been considered by other approaches as well as a family of piecewise linear 
systems known to be stable and  with the number of state variables increasing from $d = 2$ to 
$d = 9$. We show that our approach terminates faster than the conservative upper bounds established 
by our theoretical analysis. 

\subsection{Related Work}

Piecewise linear dynamical systems appear naturally in a wide range of applications, or as 
approximations of nonlinear systems~\citep{pettit1995analyzing,sun2011stability}.
Deciding stability of  such systems is known to be extremely challenging even when restricted to 
simple dynamics~\citep{prabhakar2013onthedecidability}.
Ahmadi and Jungers show that for
switched linear systems (even in $2$ dimensions and with $2$ modes), there is no
a priori bound on the complexity of polynomial, polyhedral or piecewise quadratic
functions sufficient to prove stability of the system~\citep{ahmadi2016lower}.
Moreover, the problem of deciding if there exists a polyhedral Lyapunov function with fixed number of facets
is NP-hard~\citep{berger2022learning}.
This motivates the benefit of computational methods, like the one described in this work,
that do not restrict the complexity of the  Lyapunov function up front.

Besides the many approaches for studying the stability of linear hybrid systems that
consider quadratic or polynomial Lyapunov functions~\citep{hassibi1998quadratic,sun2011stability,jungers2009thejoint},
polyhedral Lyapunov functions have also been well studied for this purpose. 
Available techniques for the computation of polyhedral Lyapunov 
functions can be divided into optimization-based and set-theoretic methods. Set-theoretic methods
typically proceed by computing the image of some polyhedral set in the state space by the system and updating 
this set until it eventually provides a polyhedral Lyapunov 
function for the system~\citep{miani2005maxisg,blanchini2015settheoretic,guglielmi2017polytope}.
These methods are however usually restricted to discrete-time systems and often lack of formal complexity bounds.
Optimization-based methods~\citep{polanski2000onabsolute,lazar2011oninfinity,ambrosino2012aconvex,kousoulidis2021polyhedral}
aim to solve the nonlinear, nonconvex optimization problem accounting for the existence of a polyhedral Lyapunov function.
This is achieved for instance by considering convex relaxations of the problem~\citep{ambrosino2012aconvex}, fixing 
some variables to make the problem  tractable~\citep{polanski2000onabsolute,lazar2011oninfinity},
or using an alternating descent procedure to solve the optimization problem by repeatedly fixing one set of variables 
to fixed values while minimizing over the remaining variables~\citep{kousoulidis2021polyhedral}.
These methods are conservative as they lack formal guarantees of finding a global optimum for the  optimization problem, 
so that they can fail to find a polyhedral Lyapunov function even if the system admits one.
Our approach borrows from both of the above ones: namely, it uses an iterative process to recursively 
update the candidate Lyapunov function, and the update is done by solving a linear programming problem.
Moreover, the approach allows to study continuous-time systems and comes with formal guarantees of convergence 
and complexity bounds on the running time of the procedure.

Our approach generates a finite set of states and uses them to infer a Lyapunov function candidate. 
This forms the basis of a ``learning'' approach to finding Lyapunov functions that has also been 
well studied. Topcu et al.~use a randomly sampled set of states to derive polynomial Lyapunov functions 
with Sum-of-Squares programming to verify the result~\citep{topcu2008local}. 
A key difference in this paper is that the generation of our witness set is based on 
counterexamples that refute previous candidates. Futhermore, we introduce a ``gap'' in 
our formulation of the learning versus the verification problems in order to yield formal 
termination guarantees. 
Kapinski et al.~learn Lyapunov functions for nonlinear systems 
by combining learning of polynomial functions from finitely many samples~\citep{kapinski2014simulationguided}.
Recently, there have been many approaches that seek to infer neural networks as Lyapunov functions 
simultaneously with neural network controllers for dynamical systems~\citep{chang2019neural,dai2021lyapunovstable}.
Likewise, Abate et al.~present an approach that synthesizes neural-network-based Lyapunov functions for nonlinear systems~\citep{abate2021formal}.
In all these approaches, a form of back-propagation is used to  learn neural networks and a suitable solver is used to 
verify the result: dReal solver for nonlinear constraints in the case of~\citet{chang2019neural}, Satisfiability-Modulo 
Theory (SMT) solver in the case of~\citet{abate2021formal} and a mixed integer solver 
in the case of~\citet{dai2021lyapunovstable}. Our approach does not synthesize controllers (although we would seek 
to do so in our future work). At the same time, neural networks with activation functions such as ReLU can be 
regarded as piecewise linear systems. 
Many of the approaches cited above use heuristic global optimization to generate counterexamples, whereas some 
approaches (notably, \citealp{chang2019neural,dai2021lyapunovstable,abate2021formal}) use expensive nonconvex solvers 
to find counterexamples. 
Nevertheless,  there are no guarantees that this approach would find a Lyapunov function, 
even if one exists. In contrast, our approach offers a  guarantee subject to 
constraints on the robustness and eccentricity of the desired function.  
Ravanbakhsh et al.~derive polynomial Lyapunov functions by using a strategy 
similar to this paper: they generate polynomial Lyapunov functions for 
nonlinear systems by iterating between learning from a finite set of 
witnesses to using Sum-of-Squares optimization for verification~\citep{ravanbakhsh2019learning}.
Unlike Ravanbakhsh et al.,\ our work focuses on linear hybrid systems and polyhedral Lyapunov 
functions. This helps us avoid the use of Sum-of-Squares relaxations in favor of linear programming problems 
that lend themselves to precise and efficient solvers. As a result, our approach provides guarantees 
upon termination that are stronger than what could be obtained by a direct application of the methods 
discussed above to piecewise linear dynamics.

This work also extends from our recent work that provides a counterexample-guided method
to compute polyhedral Lyapunov functions with fixed number of facets for
linear hybrid dynamics~\citep{berger2022learning}.
The most important  difference of the present work  is that our method
does not place any a priori bound on the complexity of the learned Lyapunov function. Furthermore, 
our method presented here provides a guarantee that if it fails then no polyhedral Lyapunov function satisfying 
user-specified limits on robustness exists for the underlying system.

\noindent\emph{Outline.}  The paper is organized as follows.
In Section~\ref{sec-problem-setting}, we introduce the problem of interest and remind important concepts related to Lyapunov analysis and polyhedral functions.
In Section~\ref{sec-description-algorithm}, we describe the algorithmic process to compute polyhedral Lyapunov functions.
In Section~\ref{sec-termination-proof}, we provide a proof of termination and soundness of the algorithm.
Finally, in Section~\ref{sec-examples-applications}, we demonstrate the applicability of the process on numerical examples.

\noindent\emph{Notation.} $\lVert\cdot\rVert$ denotes a vector norm in $\Re^d$ (typically 
the $L_1$ or $L_{\infty}$ norm), and $\Sbb:\{x\in\Re^d:\lVert x\rVert=1\}$ is the
 associated unit sphere. $\lVert\cdot\rVert_*$ denotes the dual norm of $\lVert\cdot\rVert$,
  defined by $\lVert c\rVert_*=\max\,\{c^\top x:x\in\Sbb\}$
(e.g., if $\lVert\cdot\rVert$ is the $L_1$ norm, then $\lVert\cdot\rVert_*$ is the $L_\infty$ norm).
By extension, $\lVert\cdot\rVert$ also denotes the matrix norm induced by $\lVert\cdot\rVert$ in $\Re^\dxd$,
defined by $\lVert A\rVert=\max\,\{\lVert Ax\rVert:x\in\Sbb\}$.

%%%%%%%%%%%%%%%%%%%%%%%%%%%%%%%%%%%%%%%%%%%%%%%%%%%%%%%%%%%%%%%%%%%%%%%%%%%%%%%%
\section{Problem setting}\label{sec-problem-setting}

We study continuous-time, uncertain, piecewise linear dynamical systems.

\begin{defn}\label{def-system}
A \emph{continuous-time, uncertain, piecewise linear dynamical system} is
described by a finite set of modes $Q$, wherein each mode $q\in Q$ is associated with~(i)
a closed polyhedral conic region $H_q\subseteq\Re^d$, and~(ii) a transition matrix $A_q\in\Re^\dxd$.
\end{defn}

Let $\calF:\Re^d\rightrightarrows\Re^d$ be the set-valued function defined 
by $\calF(x)=\{A_qx:q\in Q,\:x\in H_q\}$. $\calF$ completely describes the 
dynamics of  the system. Unless otherwise stated, we will refer to the 
dynamical system as ``System $\calF$''; the sets $Q$, $(H_q)_{q\in Q}$ and 
matrices $(A_q)_{q\in Q}$ being implicit in the definition of $\calF$.

% \begin{defn}[Trajectory]\label{def-trajectory}
% A function $\xi:\Re_{\geq0}\to\Re^d$ is called a \emph{trajectory} of System $\calF$
% if it satisfies $\xi'(t)\in\calF(\xi(t))$ for almost all $t\in\Re_{\geq0}$.
% \end{defn}

A function $\xi:\Re_{\geq0}\to\Re^d$ is called a \emph{trajectory} of System $\calF$
if it satisfies $\xi'(t)\in\calF(\xi(t))$ for almost all $t\in\Re_{\geq0}$.
We refer the reader to~\citet{goebel2012hybrid} for results concerning 
the existence  and uniqueness of trajectories for piecewise, uncertain,
linear dynamical systems.

\subsection{Polyhedral Lyapunov stability analysis}

We aim to study the stability of the origin for 
System $\calF$, using Lyapunov analysis.
First, let us recall the definition of the \emph{Dini derivative}.

\begin{defn}[Dini derivative]\label{def-dini-derivative}
Let $f:\Re^d\to\Re$ and $x,v\in\Re^d$.
The \emph{upper-right Dini derivative} of $f$ at $x$ in the direction of
 $v$, denoted by $D^+f(x;v)$,
is defined by $D^+f(x;v)=\limsup_{s\to0^+}\frac{f(x+sv)-f(x)}s$.
\end{defn}

\begin{defn}[Lyapunov function]\hfill
\begin{itemize}
\item A continuous function $V:\Re^d\to\Re_{\geq0}$ is a \emph{potential Lyapunov function}
if~(i) for all $x\in\Re^d$, $V(x)=0$ iff $x=0$, and~(ii) $V$ is radially unbounded\!%
\footnote{A function $V:\Re^d\to\Re$ is \emph{radially unbounded} if $\lim_{r\to\infty}\allowbreak\inf\,\{V(x):x\in X,\,\lVert x\rVert\geq r\}=\infty$.}\!.
\item A potential Lyapunov function $V:\Re^d\to\Re_{\geq0}$ is a \emph{Lyapunov function} for System $\calF$ if
for all $x\in\Re^d\setminus\{0\}$ and $v\in\calF(x)$, $D^+V(x;v)<0$.
\end{itemize}\vskip0pt
\end{defn}

It is well known that if System $\calF$ admits a Lyapunov function,
then $0$ is an asymptotically stable equilibrium point for System $\calF$
(see, e.g.,~\citealp[Theorem 2.19]{blanchini2015settheoretic}, 
or~\citealp[Theorem 4.2]{khalil2002nonlinear}). In this paper, we look 
for \emph{polyhedral} Lyapunov functions, which are defined as the pointwise maximum
of a finite set of linear functions.

\begin{defn}[Polyhedral function]\label{def-piecewise-linear-function}
A function $V:\Re^d\to\Re$ is a \emph{polyhedral function} if there is a 
finite set $\calV\subseteq\Re^d$
s.t.~for all $x\in\Re^d$, $V(x)=\max_{c\in\calV}\,c^\top x$.
\end{defn}

As a convention, given a polyhedral function $V$, we 
let $\calV$ be the set of coefficients of the various 
linear ``pieces'' defining $V$ according to 
Definition~\ref{def-piecewise-linear-function}.
Note that $\calV$ needs not be uniquely defined for a given 
function $V$. However, this will not pose an issue 
for the approach used in this paper. 
We also let $\Vmax:\max_{c\in\calV}\lVert c\rVert_*$,
and for all $x\in\Re^d$, we let $\calV(x):\{c\in\calV:V(x)=c^\top x\}$.

\begin{prop}\label{pro-piecewise-linear-lyapunov}
A polyhedral function $V$ is a Lyapunov function for System $\calF$ iff
the following conditions hold:
\begin{itemize}[leftmargin=25pt]
\item[\upshape(C1)] For all $x\in\Re^d\setminus\{0\}$, it holds that $V(x)>0$.
\item[\upshape(C2)] For all $x\in\Re^d\setminus\{0\}$, $v\in\calF(x)$ and $c\in\calV(x)$,
it holds that $c^\top v<0$.
\end{itemize}\vskip0pt
\end{prop}

\begin{pf}
The fact that (C1) is equivalent to $V$ being a potential Lyapunov function is straightforward.
Now, for a potential Lyapunov function $V$, the fact that (C2) is equivalent to $V$ being a Lyapunov function for System $\calF$
follows from the expression of the Dini derivative of $V$: for any $x\in\Re^d$ and $v\in\Re^d$,
$D^+V(x;v)=\max\,\{c^\top v:c\in\calV(x)\}$~\citep[Eq.\ 2.30]{blanchini2015settheoretic}.\qed
\end{pf}

%%%%%%%%%%%%%%%%%%%%%%%%%%%%%%%%%%%%%%%%%%%%%%%%%%%%%%%%%%%%%%%%%%%%%%%%%%%%%%%%%%%%%%%%%%%%%%%%%%%%
\subsection{Eccentricity and Robustness}

We will now recast the conditions in Proposition~\ref{pro-piecewise-linear-lyapunov}
in a form that is tractable for optimization solvers and thus enables the overall 
approach that we will develop in this paper.  To do so, we will describe parameters 
for a Lyapunov function that will be interpreted in terms of its \emph{eccentricity}
and \emph{robustness} to perturbations.

\begin{defn}[Robust Lyapunov conditions]\label{def-robust-lyapunov-conditions}
Let $\epsilon \geq 1$ be an \emph{eccentricity} parameter, and 
$\theta>0$ and $\delta>0$ be \emph{robustness} parameters. A function $V$ 
is an \emph{$(\epsilon,\theta,\delta)$-robust Lyapunov function} iff the following 
conditions hold: 
\begin{itemize}[leftmargin=25pt]
\item[\upshape(D1)] For all $x\in\Re^d$,  $V(x)\geq\frac1\epsilon\Vmax\lVert x\rVert$.
\item[\upshape(D2)] For all $x\in\Re^d$, $v\in\calF(x)$ and $c\in\calV$, 
$c^\top v\leq\frac1\theta(V(x)-c^\top x)-\delta\Vmax\lVert x\rVert$.
\end{itemize}
\end{defn}

\begin{thm}\label{thm-equivalent-piecewise-linear-lyapunov}
Let $V$ be a polyhedral function with $\Vmax>0$.
$V$ is a Lyapunov function for System $\calF$ iff
there exist $\epsilon\geq1$, $\theta>0$ and $\delta>0$ s.t. $V$ 
is an $(\epsilon,\theta,\delta)$-robust Lyapynov function. 
\end{thm}

\begin{pf}
See Appendix~\ref{app-sec-proof-thm-equivalent-piecewise-linear-lyapunov}.\qed
\end{pf}

The parameter $\epsilon$ in Definition~\ref{def-robust-lyapunov-conditions}
measures the \emph{eccentricity} of $V$, defined as
$\frac{\max_{x\in\Re^d:V(x)=1}\,\lVert x\rVert}{\min_{x\in\Re^d:V(x)=1}\,\lVert x\rVert}$.
We say that $V$ has eccentricity $\epsilon$ if the latter ratio is smaller than or equal to $\epsilon$.
Naturally, the eccentricity $\epsilon$ for any function $V$ satisfies $\epsilon\geq1$.

% The parameter $\epsilon$ in Theorem~\ref{thm-equivalent-piecewise-linear-lyapunov}
% measures the \emph{eccentricity} of $V$.

% \begin{defn}[Eccentricity]\label{def-eccentricity}
% We say that $V$ has \emph{eccentricity} (at most) $\epsilon$ if 
% $\epsilon\geq\frac{\max_{x\in\Re^d:V(x)=1}\,\lVert x\rVert}{\min_{x\in\Re^d:V(x)=1}\,\lVert x\rVert}$.
% \end{defn}

% Naturally, the eccentricity $\epsilon$ for any function $V$ satisfies $\epsilon \geq 1$.

\begin{lem}\label{lem-eccentricity}
Let $\epsilon\geq1$.
Then, $\epsilon$ satisfies (D1) in Definition~\ref{def-robust-lyapunov-conditions} iff
 $V$ has eccentricity $\epsilon$.
\end{lem}

\begin{pf}
Direct from (D1) in Definition~\ref{def-robust-lyapunov-conditions}
and the observation that $\min_{x\in\Re^d:V(x)=1}\,\lVert x\rVert=1/\Vmax$.\qed
\end{pf}

The parameter $\theta$ in Definition~\ref{def-robust-lyapunov-conditions}
can be seen as a \emph{time discretization} parameter, as (D2) can be
rewritten as $c^\top(x+\theta v)\leq V(x)-\theta\delta\Vmax\lVert x\rVert$.
When combined with the eccentricity $\epsilon$ and ``slack'' $\delta$, the
largest time step that can be used in (D2) also gives a measure of the \emph{robustness} of $V$
as a Lyapunov function w.r.t.~perturbations of System $\calF$, as we will see below.

% The parameters $(\theta,\delta)$ in Definition~\ref{def-robust-lyapunov-conditions} measure
% the \emph{robustness} of $V$ as a Lyapunov function w.r.t.~perturbations of System $\calF$.
% That is, it measures the largest perturbations that can be applied on the regions and vector fields of the system
% s.t.~$V$ remains a Lyapunov function for the perturbed system.
% The following definitions and theorems formalize this result.
% Denote $\fmax:\max\,\{\lVert v\rVert:x\in\Sbb,\:v\in\calF(x)\}$.

First, we establish the following result that holds in the absence of any perturbations. It 
notes that an $(\epsilon,\theta,\delta)$-robust Lyapunov function $V$ remains a robust 
Lyapunov function for $\epsilon'\geq\epsilon$, $\theta'\leq\theta$ and $\delta'\leq\delta$.

\begin{lem}
Let $V$ be an $(\epsilon,\theta,\delta)$-robust Lyapunov function.
It holds that $V$ is an $(\epsilon',\theta',\delta')$-robust for any 
$\epsilon'\in[\epsilon,\infty)$, $\theta'\in(0,\theta]$ and $\delta'\in(0,\delta]$.
\end{lem}

\begin{pf}
Suppose $V$ satisfies (D1) and (D2) in Definition~\ref{def-robust-lyapunov-conditions} 
for parameters $(\epsilon,\theta,\delta)$. Since $\epsilon'\geq\epsilon$, we 
obtain from (D1) that for all $x\in\Re^d$,
$V(x)\geq\frac1\epsilon\Vmax\lVert x\rVert\geq\frac{1}{\epsilon'} \Vmax\lVert x\rVert$. 
Likewise, from (D2), we obtain that for all $x\in\Re^d$, $v\in\calF(x)$ and $c\in\calV$,
$c^\top v \leq \frac1\theta(V(x)-c^\top x)-\delta\Vmax\lVert x\rVert
\leq\frac1{\theta'}(V(x)-c^\top x)-\delta\Vmax\lVert x\rVert
\leq \frac{1}{\theta'} (V(x)-c^\top x) - \delta' \Vmax\lVert x\rVert$,
where in the second inequality, we have used that $V(x)-c^\top x\geq0$ and $\theta'\leq\theta$,
and in the third inequality, we have used that $\delta'\leq\delta$.
Therefore, $V$ satisfies (D1)--(D2) in
Definition~\ref{def-robust-lyapunov-conditions} for $(\epsilon',\theta',\delta')$.\qed
\end{pf}

We now introduce the notion of \emph{perturbation} of System $\calF$, which will lead
to the notion of \emph{robustness} of a Lyapunov function w.r.t.~perturbations of the system.
For the sake of simplicity, we assume for this definition and
Theorems~\ref{thm-robustness-sufficient} and~\ref{thm-robustness-necessary} below
that all regions $(H_q)_{q\in Q}$ satisfy $H_q=\Re^d$.
% Thus, we only consider perturbations of the transition matrices $A_q$.
Let us define $\amax:\max_{q\in Q}\,\lVert A_q\rVert$.

\begin{defn}[Perturbed system]\label{def-perturbed-system}
A \emph{$\gamma$-perturbation} of System $\calF$ (with $H_q=\Re^d$) is a 
continuous-time, uncertain, piecewise linear dynamical system (Definition~\ref{def-system}) with
set of modes $Q'=Q$, set of regions $(H'_q)_{q\in Q}$ satisfying $H'_q=\Re^d$ for all $q\in Q$,
and set of matrices $(A'_q)_{q\in Q}$ satisfying $\lVert A'_q-A_q^{}\rVert\leq\gamma\amax$
for all $q\in Q$.
\end{defn}

% A Lyapunov function for System $\calF$ is \emph{robust to small perturbations of the
% system} if it is a Lyapunov function for any $\gamma$-perturbation of System $\calF$
% with $\gamma>0$ small.

The following theorems establish links between the robustness parameters $(\epsilon,\theta,\delta)$
of the polyhedral Lyapunov function and its robustness w.r.t.~$\gamma$-perturbations of the system.

\begin{thm}[Sufficient condition for robustness]\label{thm-robustness-sufficient}
Let $V$ be an $(\epsilon,\theta,\delta)$-robust polyhedral Lyapunov function for system $\calF$. 
Then, $V$ is a Lyapunov function for any $\gamma$-perturbation of System $\calF$
with $\gamma\in(0,\frac\delta\amax)$.
\end{thm}

\begin{pf}
Let $\calF'$ be a $\gamma$-perturbation of System $\calF$.
Fix $x\in\Sbb$, $v'\in\calF'(x)$ and $c\in\calV(x)$.
According to Definition~\ref{def-perturbed-system}, let $v\in\calF(x)$
be s.t.~$\lVert v-v'\rVert\leq\gamma\amax$.
It holds that $\lvert c^\top v-c^\top v'\rvert\leq\gamma\amax\Vmax$.
Hence, $c^\top v' \leq
c^\top v + \gamma\amax\Vmax <
c^\top v + \delta\Vmax$,
where the last inequality comes from the assumption on $\gamma$.
By (D2), it follows that $c^\top v'<0$,
so that (C2) in Proposition~\ref{pro-piecewise-linear-lyapunov} is satisfied for $x$ and $v'$.
Since $x$ and $v'$ were arbitrary, this concludes the proof.\qed
\end{pf}

% \begin{pf}
% See Appendix~\ref{app-sec-proof-thm-robustness-sufficient}.\qed
% \end{pf}

\begin{thm}[Necessary condition for robustness]\label{thm-robustness-necessary}
Let $\epsilon\geq1$ and $\gamma>0$.
Let $V$ be a polyhedral potential Lyapunov function with eccentricity $\epsilon$.
Assume that $V$ is a Lyapunov function for any $\gamma$-perturbation of System $\calF$.
Then, (D2) in Definition~\ref{def-robust-lyapunov-conditions} holds with any $(\theta,\delta)$
satisfying $0<\theta\leq\frac\gamma{2\amax(\epsilon+\gamma)}$ and $0<\delta\leq\frac{\gamma\amax}{2\epsilon}$.
\end{thm}

\begin{pf}
See Appendix~\ref{app-sec-proof-thm-robustness-necessary}.\qed
\end{pf}

% The general case, in which we consider perturbations of the regions $H_q$,
% is similar with one extra parameter for the ``region perturbation'' bound.

The relation between the parameters $(\epsilon,\theta,\delta)$ and the robustness of the Lyapunov function
w.r.t.~perturbations of the system being established, we focus in the following of the paper on finding
$(\epsilon,\theta,\delta)$-robust polyhedral Lyapunov functions for System $\calF$.
% A pair of robustness parameters $(\theta_1,\delta_1)$ is \emph{no worse} than another one $(\theta_2,\delta_2)$ if
% $\theta_1\geq\theta_2$ and $\delta_1\geq\delta_2$.
% We write this as $(\theta_1,\delta_1)\geq(\theta_2,\delta_2)$.
% If $(\theta_1,\delta_1)\ngeq(\theta_2,\delta_2)$, then, without further information,
% we cannot conclude whether $(\theta_1,\delta_1)$ or $(\theta_2,\delta_2)$ would give better properties
% in terms of robustness w.r.t.~perturbations of the system.

Finally, we observe that the parameters $\epsilon$, $\theta$ and $\delta$
are invariant w.r.t.~positive scaling of $V$.
That is, if $V$ is an $(\epsilon,\theta,\delta)$-robust polyhedral Lyapunov function,
then so is the function $\frac1\lambda V$ for any $\lambda>0$.
Therefore, in the following, we restrict our attention to polyhedral functions with $\Vmax\leq1$,
i.e., with $\calV\subseteq\Bbb^*:\{c\in\Re^d:\lVert c\rVert_*\leq1\}$.

\begin{exmp}[Running illustrative example]\label{exa-running}
Consider the continuous-time piecewise linear dynamical system described by
$Q:\{1,2\}$, $H_1:\Re\times\Re_{\geq0}$, $H_2:\Re\times\Re_{\leq0}$,
\[
A_1 : \left[\begin{array}{cc} -0.2 & 1.0 \\ -1.0 & -0.2 \end{array}\right] \quad\text{and}\quad
A_2 : \left[\begin{array}{cc} 0.01 & 1.0 \\ -1.0 & 0.01 \end{array}\right].
\]
The vector field of the system is represented in Figure~\ref{fig-illustrative-generator-verifier} (gray arrows).
This system does not admit a polynomial Lyapunov function; see Lemma~\ref{lem-exa-running} below.
Nevertheless, as we will see throughout the paper, we can compute a polyhedral Lyapunov function for the system,
thereby proving that it is asymptotically stable.
We can also evaluate the stability margin of the system,
by computing parameters $(\epsilon,\theta,\delta)$ for which the system does not 
admit a robust polyhedral Lyapunov function.
\end{exmp}

\begin{lem}\label{lem-exa-running}
System $\calF$ described in Example~\ref{exa-running} does not admit a polynomial Lyapunov function.
\end{lem}

\begin{pf}
For a proof by contradiction, assume that $V$ is a polynomial Lyapunov function for the system.
Let $x_0:[1,0]^\top$, and let $p:\Re\to\Re$ be the univariate polynomial defined by $p(r)=V(rx_0)$.
Since $V$ is radially unbounded, $p$ is of nonzero even degree.
Fix $r\in\Re_{>0}$ and let $x:\Re_{\geq0}\to\Re^d$ be the trajectory with $x(0)=rx_0$.
It holds that $x(\pi)=e^{A_2\pi}rx_0=-e^{0.01\pi}rx_0$.
Hence, $V(x(0))=p(r)$ and $V(x(\pi))=p(-e^{0.01\pi}r)$.
% Note that $e^{0.01 \pi}>1$.
Since $p$ is of nonzero even degree and $e^{0.01 \pi}>1$, there is $r>0$ s.t.~$p(-e^{0.01\pi}r)>p(r)$.
Thus, with such a $r$, $V(x(\pi))>V(x(0))$, contradicting the property that $V$ decreases
along the trajectories of the system~\citep[Theorem 4.1]{khalil2002nonlinear},
concluding the proof.\qed
\end{pf}

%%%%%%%%%%%%%%%%%%%%%%%%%%%%%%%%%%%%%%%%%%%%%%%%%%%%%%%%%%%%%%%%%%%%%%%%%%%%%%%%%%%%%%%%%%%%%%%%%%%%
\subsection{Overview of the Algorithm}

We are given as inputs~(i) the description of the system $\calF$, and~(ii) parameters $(\epsilon,\theta,\delta)$.
This paper will present an algorithm to~(a) compute a polyhedral Lyapunov function for System $\calF$ if one exists;
or~(b) conclude that no polyhedral Lyapunov function $V$ exists 
with eccentricity $\epsilon$ and robustness parameters $(\theta,\delta)$.

As expected, we conclude that the system is stable if the procedure computes a polyhedral Lyapunov function.
Furthermore, the parameters associated with this function helps determine how much the system can 
be perturbed while guaranteeing stability (as given by Theorem~\ref{thm-robustness-necessary}). 
However, if the procedure fails to find a Lyapunov function, we conclude 
that no Lyapunov function with eccentricity smaller than $\epsilon$ 
and robustness parameters larger than $\theta,\delta$ exists. 
Thus, we conclude upon failure of our algorithm 
that the system itself is unstable or its stability requires 
a polyhedral Lyapunov function with larger eccentricity or smaller robustness
than the input limits provided.

The iterative process is based on maintaining a finite set of states, 
called \emph{witnesses}: $X:\{x_1,\ldots,x_N\}\subseteq\Re^d$.
The witness set is initialized to $X:\emptyset$.
Each step iterates between two algorithms in succession:
\begin{itemize}
\item A \emph{Learner}, which computes a polyhedral function satisfying
the conditions of Theorem~\ref{thm-equivalent-piecewise-linear-lyapunov} over the finite set of witnesses,
or concludes that no polyhedral Lyapunov function satisfying these conditions for all $x\in\Re^d$ exists,
that is, there is no polyhedral Lyapunov function for System $\calF$
with eccentricity $\epsilon$ and robustness parameters $(\theta,\delta)$.
\item A \emph{Verifier}, which, given a candidate polyhedral function (found by the learner),
verifies whether the conditions of Proposition~\ref{pro-piecewise-linear-lyapunov} are satisfied by all states.
If the verification succeeds, we have our desired Lyapunov function.
Otherwise, the verifier algorithm returns a point (called a \emph{counterexample})
where the candidate fails to satisfy the Lyapunov conditions.
\end{itemize}

At the end of each iteration, we have three possible outcomes:~(a)
the learner refutes the existence of a ``robust-enough'' Lyapunov function,
i.e., with eccentricity $\epsilon$ and robustness parameters $(\theta,\delta)$;~(b)
the candidate polyhedral function verifies the Lyapunov conditions;
or~(c) a new witness point is added and the set $X$ grows in cardinality by one element.
We demonstrate that this process will eventually terminate,
and provide upper bounds on the total number of iterations to termination.

%%\comgb{Can we improve the paragraph below?}

% {\new

% Let us stress out the gap in the output of the algorithm,
% which either refutes the existence of a Lyapunov functions with parameters $\epsilon$, $\theta$ and $\delta$,
% or outputs a polyhedral function that is a Lyapunov function for System $\calF$, but not necessarily with parameters $\epsilon$, $\theta$ and $\delta$.
% The algorithm can be refined to output a Lyapunov function for System $\calF$ with parameters $\epsilon_o$, $\theta_o$ and $\delta_o$, where $\epsilon_o>\epsilon$, $\theta_o<\theta$ and $\delta_o<\delta$.
% The larger the gap between the above quantities, the more efficient the algorithm.
% For the sake of simplicity, we will focus in the rest of the paper on the case where $\epsilon_o=\infty$ and $\theta_o=\delta_o=0$.

% }

The algorithm is described in Section~\ref{sec-description-algorithm}
and the proof of its termination and soundness is presented in Section~\ref{sec-termination-proof}.

%%%%%%%%%%%%%%%%%%%%%%%%%%%%%%%%%%%%%%%%%%%%%%%%%%%%%%%%%%%%%%%%%%%%%%%%%%%%%%%%
\section{Description of the algorithm}\label{sec-description-algorithm}

First, we present the learner, then the verifier, and finally the overall algorithmic process.
We assume input parameters $\epsilon\geq1$, $\theta>0$ and $\delta>0$.

%%%%%%%%%%%%%%%%%%%%%%%%%%%%%%%%%%%%%%%%%%%%%%%%%%%%%%%%%%%%%%%%%%%%%%%%%%%%%%%%
\subsection{Learner: Computation of a candidate polyhedral Lyapunov function}

Let $X:\{x_1,\ldots,x_N\}\subseteq\Re^d$ be a finite set of witnesses.
We consider the following computational problem, which aims to find a polyhedral function
satisfying the conditions of Definition~\ref{def-robust-lyapunov-conditions} at every witness point.
The function $V$ will have as many as $\lvert X\rvert$ pieces, one piece 
$c_x$ associated with each witness $x\in X$.

\begin{prob}\label{prob-learner}
Find a polyhedral function $V$ with set $\calV:\{c_x:x\in X\}\subseteq\Bbb^*$
s.t.~for all $x\in X$,~(i) $c_x^\top x\geq\frac1\epsilon\lVert x\rVert$,
and~(ii) for all $v\in\calF(x)$ and $c\in\calV$, it holds that
$c^\top v\leq\frac1\theta(c_x^\top x-c^\top x)-\delta\lVert x\rVert$.
\end{prob}

The constraints (i)--(ii) in Problem~\ref{prob-learner} are linear in the decision variables $\{c_x:x\in X\}$
and the constraint $\calV\subseteq\Bbb^*$ is convex.
Thus, Problem~\ref{prob-learner} can be cast as a convex optimization problem.
These problems can be solved efficiently and accurately using for instance interior-point algorithms~\citep{boyd2004convex}.

We have the following co-soundness result:

\begin{lem}\label{lem-learner-feasibility}
If System $\calF$ admits an $(\epsilon,\theta,\delta)$-robust polyhedral
Lyapunov function, then Problem~\ref{prob-learner} is feasible.
\end{lem}

\begin{pf}
Straightforward: it suffices to take, for each $x\in X$, $c_x\in\calV(x)$,
where $V$ is an $(\epsilon,\theta,\delta)$-robust polyhedral
Lyapunov function for System $\calF$.\qed
\end{pf}

Conversely,  any solution to Problem~\ref{prob-learner} is potentially a Lyapunov function for System $\calF$.
\emph{However, it would need be verified for all state $x\in\Re^d$, rather than just at the finite set of witnesses.}

Note that the constraints (i)--(ii) in Problem~\ref{prob-learner} are invariant w.r.t.~positive scaling of the witnesses.
Therefore, in the following, we assume w.l.o.g.~that all witnesses are normalized to the unit sphere.

\begin{assum}\label{ass-norm-x}
$X\subseteq\Sbb$.
\end{assum}

% We also briefly remark on the correspondence between the witnesses in $X$ and the vectors in $\calV$.
% Implicitly, the learner seeks to find a polyhedral function $V$ corresponding to a witness set $X$
% s.t.~every witness $x$ is associated with a piece $c_x\in\calV$.
% {\new The requirement that $V(x)=c_x^\top x$ needs not be enforced, as it will be automatically satisfied if needed,
% because it can only improve the feasibility of (i) and (ii) in Problem~\ref{prob-learner}.}

\begin{figure}
\centering
\includegraphics[width=0.75\linewidth]{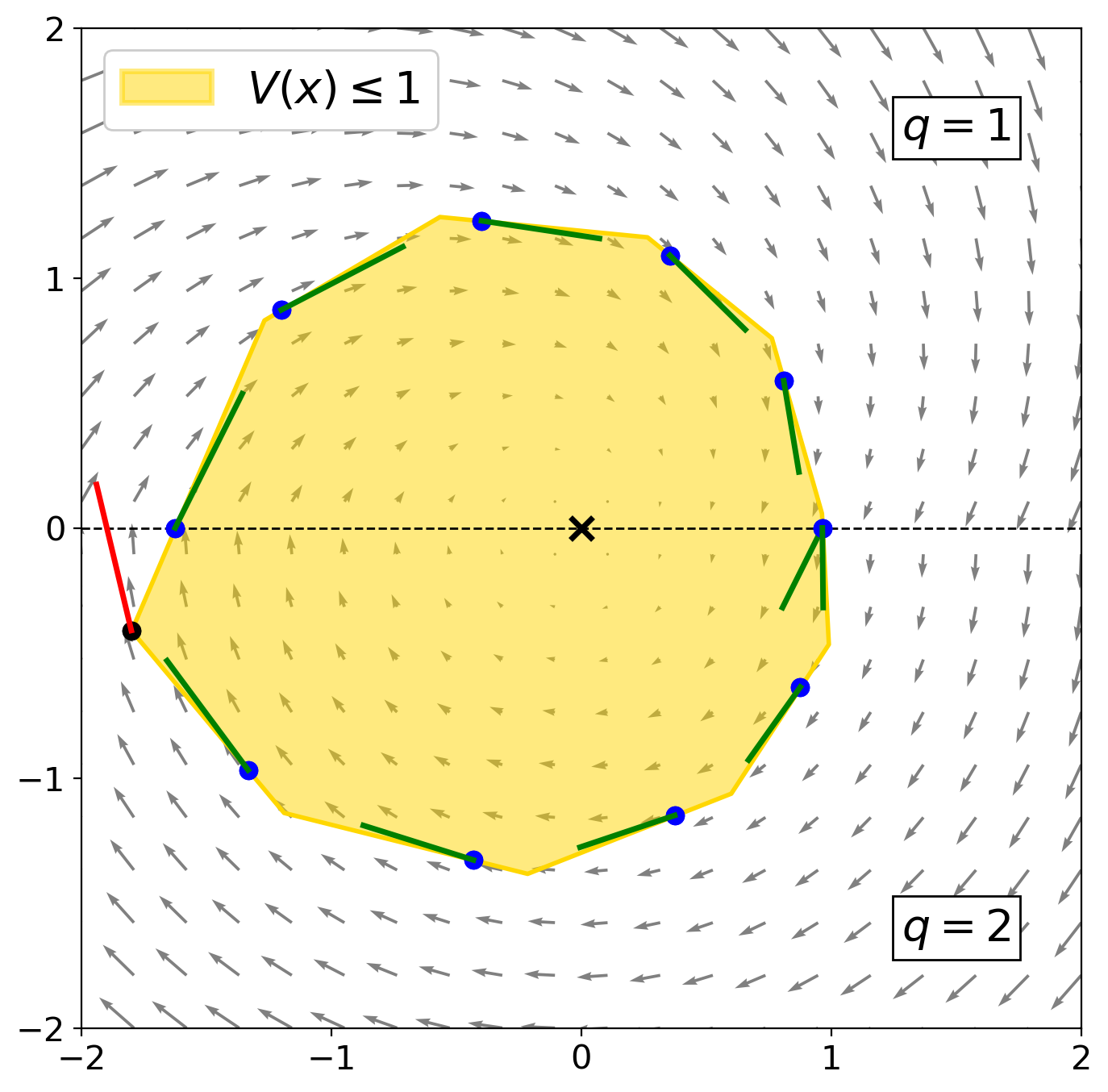}
\caption{(Running illustrative example).
Candidate polyhedral function $V$ provided by the learner (Problem~\ref{prob-learner}),
with parameters $\epsilon:10$, $\theta:1/4$ and $\delta:0.05$,
for System $\calF$ (gray arrows) and witness set $X$ (blue dots).
A counterexample (black dot) provided by the verifier (Problem~\ref{prob-verifier})
for System $\calF$ and the candidate polyhedral function $V$.}
\label{fig-illustrative-generator-verifier}
\end{figure}

\begin{exmp}[Running illustrative example]\label{exa-running-learner}
Consider System $\calF$ described in Example~\ref{exa-running}.
Let $X$ be a set of $10$ circularly placed points in $\Re^2$.
Let $\epsilon:10$, $\theta:1/4$ and $\delta:0.05$.
With these parameters, Problem~\ref{prob-learner} has a feasible solution $V$.
The $1$-sublevel set of $V$, the witness points $x\in X$ and the flow directions $v\in\calF(x)$
are represented in Figure~\ref{fig-illustrative-generator-verifier} (blue dots and green lines).
We see that the flow directions always point toward to interior of the $1$-sublevel set.
\end{exmp}

%%%%%%%%%%%%%%%%%%%%%%%%%%%%%%%%%%%%%%%%%%%%%%%%%%%%%%%%%%%%%%%%%%%%%%%%%%%%%%%%%%%%%%%%%%%%%%%%%%%%
\subsection{Verifier: Verification of the candidate polyhedral Lyapunov function and counterexample generation}

Let $V$ be a candidate polyhedral function output by the learner.
We consider the problem of computing a \emph{counterexample} $x$
that fails to satisfy the conditions (C1)--(C2) in Proposition~\ref{pro-piecewise-linear-lyapunov}.

\begin{prob}\label{prob-verifier}
Find $x\in\Re^d\setminus\{0\}$ s.t.~either~(i) $V(x)\leq0$,
or~(ii) there is $v\in\calF(x)$ and $c\in\calV(x)$ s.t.~$c^\top v\geq0$.
\end{prob}

Problem~\ref{prob-verifier} can be approached by solving the following optimization problems:~(i)
{\allowdisplaybreaks\begin{subequations}\label{eq-optim-verifier-pos}
\begin{align}
\text{find} \quad & x \\
\text{s.t.}\quad & c^\top x\leq0, \quad\forall\,c\in\calV, \label{eq-optim-verifier-pos-max} \\
& \lVert x\rVert\geq1; \label{eq-optim-verifier-pos-norm}
\end{align}
\end{subequations}}%
and~(ii) for each $c\in\calV$ and $q\in Q$:
{\allowdisplaybreaks\begin{subequations}\label{eq-optim-verifier-lie}
\begin{align}
\text{find} \quad & x\in H_q \\
\text{s.t.}\quad & c^\top x=1\geq b^\top x, \quad\forall\,b\in\calV, \label{eq-optim-verifier-lie-max} \\
& c^\top\!A_qx\geq 0. \label{eq-optim-verifier-lie-deriv}
\end{align}
\end{subequations}}

The constraint~\eqref{eq-optim-verifier-pos-norm} is nonconvex,
but it can be enforced as the disjunction of $2d$ linear constraints 
(e.g., one for each component of $x$ to be equal to $-1$ or $1$),
thereby allowing to solve~\eqref{eq-optim-verifier-pos} by solving $2d$ linear programs.
Thus, solving problems~\eqref{eq-optim-verifier-pos} and~\eqref{eq-optim-verifier-lie} amounts
to solve $2d+\lvert\calV\rvert\lvert Q\rvert$ linear programs with $d$ variables and
$\lvert\calV\rvert+\max_{q\in Q}m_q+1$ constraints
(where $m_q$ is the number of linear constraints describing $H_q$),
which can be achieved efficiently and reliably~\citep{nesterov1994interiorpoint}.

\begin{thm}[Soundness and completeness]\label{thm-verifier-complete-sound}
Problem~\ref{prob-verifier} is feasible
iff~\eqref{eq-optim-verifier-pos} or~\eqref{eq-optim-verifier-lie} has a feasible solution.
\end{thm}

\begin{pf}
Straightforward by the positive homogeneity of the problem.\qed
\end{pf}

\begin{exmp}[Running illustrative example]\label{exa-running-verifier}
Consider System $\calF$ described in Example~\ref{exa-running}.
Let $V$ be the candidate polyhedral function computed by the learner in Example~\ref{exa-running-learner}.
We check whether there is a counterexample for $V$.
It turns out to be the case since Problem~\ref{prob-verifier} has a feasible solution $x\approx[-0.62,-0.85]^\top$.
The counterexample point $x$ and the flow direction $v\in\calF(x)$
are represented in Figure~\ref{fig-illustrative-generator-verifier} (black dot and red line).
We see that the flow direction points towards the exterior of the $1$-sublevel set of $V$.
\end{exmp}

%%%%%%%%%%%%%%%%%%%%%%%%%%%%%%%%%%%%%%%%%%%%%%%%%%%%%%%%%%%%%%%%%%%%%%%%%%%%%%%%%%%%%%%%%%%%%%%%%%%%
\subsection{Overall algorithmic process}

We now describe the overall algorithmic process to compute a polyhedral Lyapunov function for System $\calF$,
or conclude that no polyhedral Lyapunov function with eccentricity $\epsilon$ and robustness parameters $(\theta,\delta)$ exists.

The process starts with an empty set of witnesses (or alternatively a finite set of witnesses given as input of the algorithm).
Then, the process enters a loop, in which at each iteration,
the following learning step and verification step are performed sequentially:~(a)
From the current set of witnesses, the learner outputs a candidate polyhedral Lyapunov functions,
or concludes that no $(\epsilon,\theta,\delta)$-robust polyhedral Lyapunov function exists.
In the latter case, the algorithm stops and outputs $\textsc{fail}$.~(b)
The verifier checks whether the candidate polyhedral function provides a valid Lyapunov function for System $\calF$.
If it is the case, then the algorithm stops and outputs the candidate polyhedral function.
Otherwise, it produces a counterexample, which is added to the witness set.
The algorithm then proceeds with the next iteration of the loop.
The process is described in Algorithm~\ref{algo-process}.

\begin{algorithm2e}[t]
\SetKw{Break}{break} 
\DontPrintSemicolon
\caption{Learning a Polyhedral Lyapunov Function.}
\label{algo-process}
\KwData{System $\calF$, eccentricity $\epsilon\geq1$, robustness parameters $\theta>0$ and $\delta>0$.}
\KwResult{Polyhedral Lyapunov function $V$, or \textsc{fail}.}
\tcc{Initialization}
$X_0\gets\emptyset$\;
\tcc{Learning loop}
\For{$k=0,1,\ldots$}{
	Compute candidate polyhedral function $V_k$ by solving Problem~\ref{prob-learner} with $X:X_k$\;
	\lIf{Problem~\ref{prob-learner} has no solution}{\Return{\textsc{fail}}}
	Find a counterexample $x_k$ for the candidate $V_k$ by solving Problem~\ref{prob-verifier} with $V:V_k$\;
	\lIf{Problem~\ref{prob-verifier} has no solution}{\Return{$V_k$}}
	$\xbar_k\gets x_k/\lVert x_k\rVert$\;
	$X_{k+1}\gets X_k\cup\{\xbar_k\}$\;
}
\end{algorithm2e}

\begin{figure}
\centering
\includegraphics[width=\linewidth]{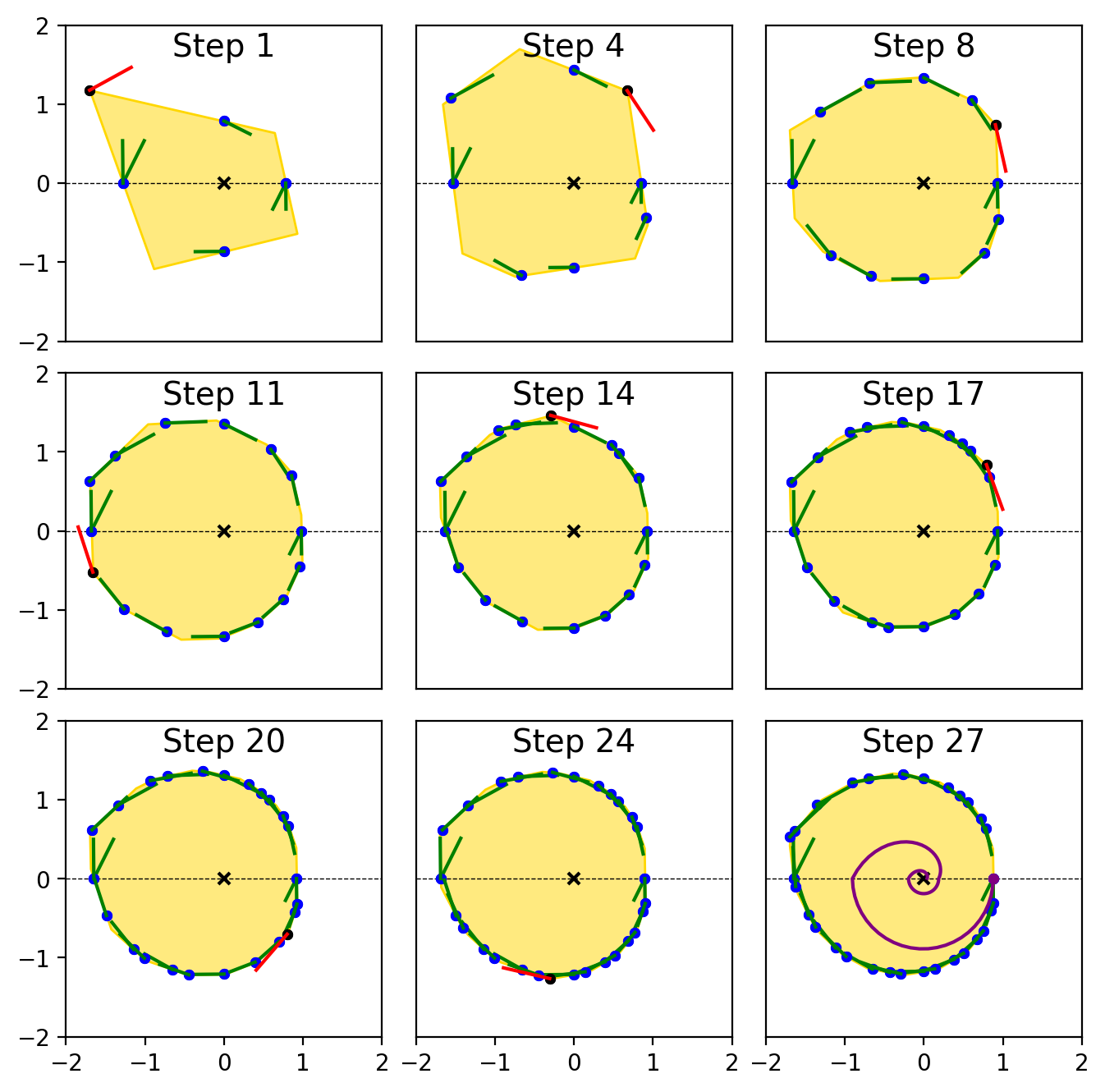}
\caption{(Running illustrative example).
Different steps of the construction of a polyhedral Lyapunov function
for System $\calF$ in Example~\ref{exa-running-process}.
\emph{Yellow:} $1$-sublevel set of $V_k$.
\emph{Blue dots:} witness points $x\in X_k$ scaled to be on the $1$-level set of $V_k$.
\emph{Green lines:} flow directions of the system at the scaled witness points.
\emph{Black dot:} counterexample points $x\in X_k$ scaled to be on the $1$-level set of $V_k$.
\emph{Red lines:} flow directions of the system at the scaled counterexample points.
\emph{Purple curve:} a sample trajectory of System $\calF$.}
\label{fig-illustrative-process}
\end{figure}

\begin{exmp}[Running illustrative example]\label{exa-running-process}
Consider System $\calF$ described in Example~\ref{exa-running}.
Let $\epsilon:10$, $\theta:1/4$ and $\delta:0.05$.
The construction of a polyhedral Lyapunov function for this system,
using Algorithm~\ref{algo-process}, is illustrated in Figure~\ref{fig-illustrative-process}.
We start with a set $X_0$ of $4$ points in $\Re^2$.
At each step $k$, a polyhedral function $V_k$,
satisfying Problem~\ref{prob-learner} with $X:X_k$, is computed.
Then, the verifier checks whether it can find a counterexample $x_k$
satisfying Problem~\ref{prob-verifier} for $V:V_k$.
For several different steps $k$, we have represented, in Figure~\ref{fig-illustrative-process},
the set $X_k$, the function $V_k$ and the counterexample $x_k$.
% More precisely, the $1$-sublevel set of $V_k$ is represented in yellow,
% the points in $X_k$, scaled to be on the boundary of the sublevel set, are represented by blue dots,
% and the flow directions of the system for each of these points are represented by green lines.
% The counterexample points $x_k$, scaled to be on the boundary of the sublevel set, are represented by black dots,
% and the flow directions of the system from these points are represented by red lines.

After 27 steps, the process has computed a polyhedral Lyapunov function for System $\calF$.
The computed function is represented in the last plot of Figure~\ref{fig-illustrative-process}.
We have also represented a sample trajectory of the system.
We observe that the trajectory does not leave the sublevel set, as predicted by the theory of Lyapunov.

Finally, if we set $\delta:0.1$ (instead of $0.05$) with the remaining 
parameters retaining their previous values ($\epsilon:10$, $\theta:1/4$),
then after $20$~steps, we conclude that  Problem~\ref{prob-learner} has no solution. Thus, the 
algorithm outputs \textsc{fail}. This shows that System $\calF$ does not admit a polyhedral Lyapunov function
with eccentricity $\epsilon \leq 10$ and robustness parameters $\theta\geq1/4$ and $\delta\geq0.1$.
\end{exmp}

%%%%%%%%%%%%%%%%%%%%%%%%%%%%%%%%%%%%%%%%%%%%%%%%%%%%%%%%%%%%%%%%%%%%%%%%%%%%%%%%
\section{Proof of termination of the algorithm}\label{sec-termination-proof}

The proof of termination exploits the gap between the constraints enforced by the learner
and the constraints verified by the verifier.
Indeed, the learner requires that the polyhedral function has eccentricity $\epsilon$ and
robustness parameters $(\theta,\delta)$ (over the set of witnesses), whereas the verifier only checks
whether the function is a valid Lyapunov function.
From this, we get that if $V$ is a candidate provided by the learner with witness set $X$,
and $X'$ is ``sufficiently close'' to the witness set $X$, then $V$ satisfies the 
conditions of the verifier on the set $X'$.
Similarly, if $V'$ is a small perturbation of $V$, then $V'$ also satisfies the conditions of the verifier
on the set $X$. This leads to two different arguments for termination with corresponding bounds on the 
number of iterations needed by our algorithm in the worst case. 
We discuss these properties in detail in Subsections~\ref{ssec-perturb-witness} and~\ref{ssec-perturb-candidate},
and derive associated complexity bounds for the algorithm.

%%%%%%%%%%%%%%%%%%%%%%%%%%%%%%%%%%%%%%%%%%%%%%%%%%%%%%%%%%%%%%%%%%%%%%%%%%%%%%%%
\subsection{Perturbation of the witness set}\label{ssec-perturb-witness}

We discuss the property that the candidate function provided by the learner
satisfies the conditions of the verifier at points that are close to the witness set.
First, we introduce the notion of ``inflation'' of the witness set.

\begin{defn}\label{def-close-witness}
Let $X\subseteq X'\subseteq\Sbb$ and $r\geq 0$.
We say that $X'$ is an \emph{$r$-inflation} of $X$ (w.r.t.~the regions $(H_q)_{q\in Q}$) if
for all $q\in Q$ and $x'\in X'\cap H_q$, there is $x\in X\cap H_q$ s.t.~$\lVert x-x'\rVert<r$.
\end{defn}

We now prove the property.
Fix $\epsilon\geq1$, $\theta>0$ and $\delta>0$.
Define $\rbar:\min\,\{\frac1\epsilon,\frac{\theta\delta}{2+\theta\amax}\}$.

\begin{thm}\label{thm-robustness-X}
Let $X\subseteq\Sbb$.
Let $V$ be a solution to Problem~\ref{prob-learner}.
Let $X'\subseteq\Sbb$ be an $\rbar$-inflation of $X$.
Then, for all $x'\in X'$,~(i) $V(x')>0$, and~(ii) for all $v'\in\calF(x')$ and $c'\in\calV(x')$,
$c'^\top v'<0$.
\end{thm}

\begin{pf}
Let $q\in Q$ and $x'\in X'\cap H_q$.
Let $x\in X\cap H_q$ be s.t.~$\lVert x-x'\rVert<\rbar$.
We first show that $V(x')>0$.
It holds that $c_x^\top x\geq\frac1\epsilon$ and
$\lvert c_x^\top x-c_x^\top x'\rvert<\rbar$.
Hence, by definition of $\rbar$, $c_x^\top x'>0$, so that $V(x')>0$.

Now, let $c\in\calV(x')$.
Let $v:A_qx$ and $v':A_qx'$.
Note that $\lVert v-v'\rVert=\lVert A_q(x-x')\rVert<\amax\rbar$.
We show that $c^\top v'<0$.
It holds that $c^\top x'\geq c_x^\top x'>c_x^\top x-\rbar$, $\lvert c^\top x-c^\top x'\rvert<\rbar$,
and $\lvert c^\top v-c^\top v'\rvert<\amax\rbar$.
Also, $c^\top v\leq\frac1\theta(c_x^\top x-c^\top x)-\delta$.
Hence, $c^\top v' < c^\top v - \frac1\theta(c_x^\top x-c^\top x) + \rbar(\frac2\theta+\amax) \leq 0$,
where the last inequality follows from the definition of $\rbar$.
Since $q$, $x$, and $c$ were arbitrary, this concludes the proof.\qed
\end{pf}

\begin{cor}\label{cor-separated-X}
Let $X_k\subseteq X_{k+1}\subseteq\Sbb$ be two consecutive witness sets
generated during the execution of Algorithm~\ref{algo-process}.
Then, $X_{k+1}$ is not an $\rbar$-inflation of $X_k$.
\end{cor}

\begin{pf}
If $X_{k+1}$ is an $\rbar$-inflation of $X_k$,
then by Theorem~\ref{thm-robustness-X}, the counterexample $\xbar_k$ does not belong to $X_{k+1}$,
contradicting the definition $X_{k+1}:X_k\cup\{\xbar_k\}$.\qed
\end{pf}

From Corollary~\ref{cor-separated-X}, we derive the following upper bound
on the number of iterations of the algorithm.
For given $s>0$ and compact set $C\subseteq\Re^d$, 
let $\mathrm{Pack}(s;C)$ denote the \emph{$s$-packing number} of $C$, i.e, 
the largest cardinality of a subset $\Chat\subseteq C$
s.t.~for all $x,y\in\Chat$, if $x\neq y$ then $\lVert x-y\rVert\geq s$.

\begin{thm}[Termination]\label{thm-termination-process-X}
Algorithm~\ref{algo-process} terminates in at most $\lvert Q\rvert\,\mathrm{Pack}(\rbar;\Sbb)$
steps.
\end{thm}

\begin{pf}
Assume that Algorithm~\ref{algo-process} produces at least $K+1$ counterexamples: $\xbar_0,\ldots,\xbar_K$.
For each $k$, let $q_k\in Q$ be s.t.~$\xbar_k\in H_{q_k}$ and $\min_{x\in X_k\cap H_{q_k}}\,\lVert x-\xbar_k\rVert\geq\rbar$
(by Corollary~\ref{cor-separated-X}, such a $q_k$ always exists).
For each $q\in Q$, let $X_K\!\downarrow\!q:\{\xbar_k:q_k=q\}$.
By the pigeonhole principle, there is $q\in Q$ s.t.~$\lvert X_K\!\downarrow\!q\rvert\geq (K+1)/Q$.
Fix such a $q$.
It holds that for all $x,y\in X_K\!\downarrow\!q$, if $x\neq y$ then $\lVert x-y\rVert\geq\rbar$.
Thus, $\lvert X_K\!\downarrow\!q\rvert$ is upper bounded by $\mathrm{Pack}(\rbar;\Sbb)$.
This proves that $K+1\leq\lvert Q\rvert\,\mathrm{Pack}(\rbar;\Sbb)$.\qed
\end{pf}

%%%%%%%%%%%%%%%%%%%%%%%%%%%%%%%%%%%%%%%%%%%%%%%%%%%%%%%%%%%%%%%%%%%%%%%%%%%%%%%%
\subsection{Perturbation of the candidate function}\label{ssec-perturb-candidate}

We discuss the property that a small enough perturbation of the candidate function
provided by the learner, still satisfies the conditions of the verifier on the same witness set.
We first formalize the property of ``small enough'' perturbation.

\begin{defn}\label{def-close-function}
Let $V_1$ and $V_2$ be two polyhedral functions and $s\geq0$.
We say that $V_1$ and $V_2$ are \emph{$s$-close} if
for all $c_1\in\calV_1$, $\inf_{c_2\in\calV_2}\,\lVert c_1-c_2\rVert_*<s$,
and for all $c_2\in\calV_2$, $\inf_{c_1\in\calV_1}\,\lVert c_1-c_2\rVert_*<s$,
\end{defn}

We now prove the property.
Fix $\epsilon\geq1$, $\theta>0$ and $\delta>0$.
Remember the definition $\rbar:\min\,\{\frac1\epsilon,\frac{\theta\delta}{2+\theta\amax}\}$.

\begin{thm}\label{thm-robustness-V}
Let $X\subseteq\Sbb$.
Let $V$ be a solution to Problem~\ref{prob-learner}.
Let $V'$ be a polyhedral function s.t.~$V$ and $V'$ are $\rbar$-close.
Then, for all $x\in X$,~(i) $V'(x)>0$, and~(ii) for all $v\in\calF(x)$ and $c'\in\calV'(x)$,
$c'^\top v<0$.
\end{thm}

\begin{pf}
Let $x\in X$.
We first show that $V'(x)>0$.
Let $c\in\calV(x)$, and let $c'\in\calV'$ be s.t.~$\lVert c-c'\rVert_*<\rbar$.
It holds that $c^\top x\geq\frac1\epsilon$ and
$\lvert c^\top x-c'^\top x\rvert<\rbar$.
Hence, by definition of $\rbar$, $c'^\top x>0$, so that $V'(x)>0$.

Now, let $v\in\calF(x)$ and $c'\in\calV'(x)$.
We show that $c'^\top v<0$.
Let $c\in\calV$ be s.t.~$\lVert c-c'\rVert_*<\rbar$,
and $c_x'\in\calV'$ be s.t.~$\lVert c_x^{}-c_x'\rVert_*<\rbar$.
It holds that $c'^\top x\geq c_x'^\top x>c_x^\top x-\rbar$, $\lvert c^\top x-c'^\top x\rvert<\rbar$ and
$\lvert c^\top v-c'^\top v\rvert<\rbar\amax$.
Also, $c^\top v\leq\frac1\theta(c_x^\top x-c^\top x)-\delta$.
Hence, $c'^\top v < c^\top v - \frac1\theta(c_x^\top x-c^\top x) + \rbar(\frac2\theta+\amax) \leq 0$,
where the last inequality follows from the definition of $\rbar$.
Since $x$, $v$ and $c'$ were arbitrary, this concludes the proof.\qed
\end{pf}

\begin{cor}\label{cor-separated-V}
Let $V_0,V_1,\ldots$ be the candidate polyhedral functions generated by the learner
during the execution of Algorithm~\ref{algo-process}.
Then, for any $k_1<k_2$, $V_{k_1}$ and $V_{k_2}$ are not $\rbar$-close.
\end{cor}

\begin{pf}
If $V_{k_1}$ and $V_{k_2}$ are $\rbar$-close,
then by Theorem~\ref{thm-robustness-V} the counterexample $\xbar_{k_1}$ is not in $X_{k_2}$,
contradicting the definition $X_{k_2}:X_{k_1}\cup\{\xbar_{k_1},\ldots,\xbar_{k_2-1}\}$.\qed
\end{pf}

From Corollary~\ref{cor-separated-V}, we derive the following upper bound
on the complexity of the algorithm. Let us define the \emph{covering number}
$\mathrm{Cov}_*(s;C)$ for a compact set $C \subseteq\Re^d$ 
as the smallest cardinality of a set $\Chat\subseteq\Re^d$
s.t.~for all $c\in\Chat$ there exists $c'\in\Chat$ s.t.~$\lVert c-c'\rVert\leq s$.

\begin{thm}[Termination]\label{thm-termination-process-V}
Algorithm~\ref{algo-process} terminates in at most $N$ steps,
where $N$ is the maximal cardinality of a set $\{V_k\}_{k=1}^N$ of polyhedral functions
satisfying Problem~\ref{prob-learner} (for some witness sets $\{X_k\}_{k=1}^N$),
that are not $\rbar$-close to each other.
It holds that $N\leq2^{\mathrm{Cov}_*(\frac\rbar2;\Bbb^*)}$.
\end{thm}

\begin{pf}
The bound $N$ on the number of steps follows directly from
Corollary~\ref{cor-separated-V}.
We derive the bound $N\leq2^{\mathrm{Cov}_*(\frac\rbar2;\Bbb^*)}$ as follows.
Let $\{V_k\}_{k=1}^N$ be as in the statement of the theorem,
and let $\Chat\subseteq\Re^d$ be an $\frac\rbar2$-covering of $\Bbb^*$.
For each $k$, let $\Chat_k\subseteq\Chat$ be a minimal $\frac\rbar2$-covering of $\calV_k$.
We show that for any $k_1\neq k_2$, $\Chat_{k_1}\neq\Chat_{k_2}$.
For a proof by contradiction, let $k_1\neq k_2$ be s.t.~$\Chat_{k_1}=\Chat_{k_2}$.
Let $c_1\in\calV_1$.
There is $c\in\Chat_{k_1}=\Chat_{k_2}$ s.t.~$\lVert c-c_1\rVert_*<\frac\rbar2$.
Moreover, since $\Chat_{k_2}$ is minimal,
there is $c_2\in\calV_2$ s.t.~$\lVert c-c_2\rVert_*<\frac\rbar2$.
For such a $c_2$, it holds that $\lVert c_1-c_2\rVert_*<\rbar$.
Similarly, one can show that for all $c_2\in\calV_2$, there is $c_1\in\calV_1$ s.t.~$\lVert c_1-c_2\rVert_*<\rbar$.
Thus, $V_1$ and $V_2$ are $\rbar$-close, a contradiction.
Hence, we have shown that $k_1\neq k_2$ implies $\Chat_{k_1}\neq C_{k_2}$.
Since there are at most $2^{\lvert\Chat\rvert}$ subsets of $\Chat$,
it follows that $K\leq 2^{\lvert\Chat\rvert}$.
Since $\Chat$ was arbitrary, this concludes the proof.\qed
\end{pf}

% Now, the fact that the algorithm outputs \textsc{fail} only if System $\calF$
% does not admit a polyhedral Lyapunov function with eccentricity $\epsilon$ and robustness parameters $(\theta,\delta)$
% follows from Theorem~\ref{thm-learner-termination}.
% Otherwise, if the algorithm provides a polyhedral function, it means that this function has passed the verifier test,
% so that it is a valid Lyapunov function for System $\calF$.\qed

Thus, we have provided two proofs of terminaton for Algorithm~\ref{algo-process} 
that yield two  bounds on the running time of the algorithm.

\begin{cor}[Termination]\label{cor-termination-process}
Algorithm~\ref{algo-process} terminates in at most
$\min\,\{\lvert Q\rvert\,\mathrm{Pack}(\rbar;\Sbb),N\}$ steps,
where $N$ is defined in Theorem~\ref{thm-termination-process-V}.

Moreover, upon termination, the process outputs \textsc{fail} only if
System $\calF$ does not admit an $(\epsilon,\theta,\delta)$-robust polyhedral Lyapunov function;
otherwise, it outputs a polyhedral Lyapunov function for System $\calF$. 
\end{cor}

\begin{pf}
The bound on the number of steps follows from
Theorems~\ref{thm-termination-process-X} and~\ref{thm-termination-process-V}.
Now, the fact that the algorithm outputs \textsc{fail} only if System $\calF$
does not admit an $(\epsilon,\theta,\delta)$-robust polyhedral Lyapunov function
follows from Lemma~\ref{lem-learner-feasibility}.
Otherwise, if the algorithm provides a polyhedral function,
it means that this function has passed the verifier test,
so that it is a valid Lyapunov function for System $\calF$ (Theorem~\ref{thm-verifier-complete-sound}).\qed
\end{pf}

\begin{rem}
For reasonable numbers of modes $\lvert Q\rvert$,
the upper bound $2^{\mathrm{Cov}_*(\frac\rbar2;\Bbb^*)}$ on $N$
is always worse than the bound $\lvert Q\rvert\,\mathrm{Pack}(\rbar;\Sbb)$.
However, it should be noted that $N$ is in general much smaller than $2^{\mathrm{Cov}_*(\frac\rbar2;\Bbb^*)}$
because the set of polyhedral functions satisfying Problem~\ref{prob-learner}
is much smaller than the set of all polyhedral functions;
especially if the system is not robustly stable.
Therefore, we have kept both bounds in Corollary~\ref{cor-termination-process}:
the first bound is more useful when the system is robustly stable, because in that case
$\rbar$ is expected to be large, so that $\mathrm{Pack}(\rbar;\Sbb)$ is smaller;
the second bound is more useful when the system is not robustly stable, because in that case
the set of possible candidate functions is expected to be smaller, so that $N$ is smaller.
\end{rem}

%%%%%%%%%%%%%%%%%%%%%%%%%%%%%%%%%%%%%%%%%%%%%%%%%%%%%%%%%%%%%%%%%%%%%%%%%%%%%%%%%%%%%%%%%%%%%%%%%%%%
\section{Numerical examples}\label{sec-examples-applications}

All computations were made on a laptop with processor Intel Core i7-7600u and 16 GB RAM running Windows.
We use Gurobi\textsuperscript{TM}, under academic license, as linear optimization solver.

%%%%%%%%%%%%%%%%%%%%%%%%%%%%%%%%%%%%%%%%%%%%%%%%%%%%%%%%%%%%%%%%%%%%%%%%%%%%%%%%%%%%%%%%%%%%%%%%%%%%
\subsection{Benchmark: $2D$ uncertain linear system}\label{ssec-benchmark-zelentsowsky}

This system was introduced in~\citet{zelentsovsky1994nonquadratic} and was used by many
authors as a benchmark comparing the performance of various Lyapunov synthesis approaches.
The system is described by $\calF(x) = \{A_px:p\in\{0,\alpha\}\}$ where
\[
A_p : \left[\begin{array}{cc}
0 & 1 \\ -2 & -1
\end{array}\right] +
p\left[\begin{array}{cc}
0 & 0 \\ -1 & 0
\end{array}\right].
\]
\citet{zelentsovsky1994nonquadratic} shows that the system with $\alpha:3.82$ admits a quadratic Lyapunov function;%
~\citet{blanchini1996onthetransient} provide a polyhedral Lyapunov function for the system with $\alpha:6$;%
~\citet{xie1997piecewise} provide a piecewise quadratic function for the system with $\alpha:6.2$;%
~\citet{chesi2009homogeneous} provide a polynomial Lyapunov function of degree $20$ for the system with $\alpha:6.8649$;
and~\citet{ambrosino2012aconvex} provide a polyhedral Lyapunov function with $9694$ vertices for the system with $\alpha:6.87$.

\begin{figure}
\centering
\includegraphics[width=0.6\linewidth]{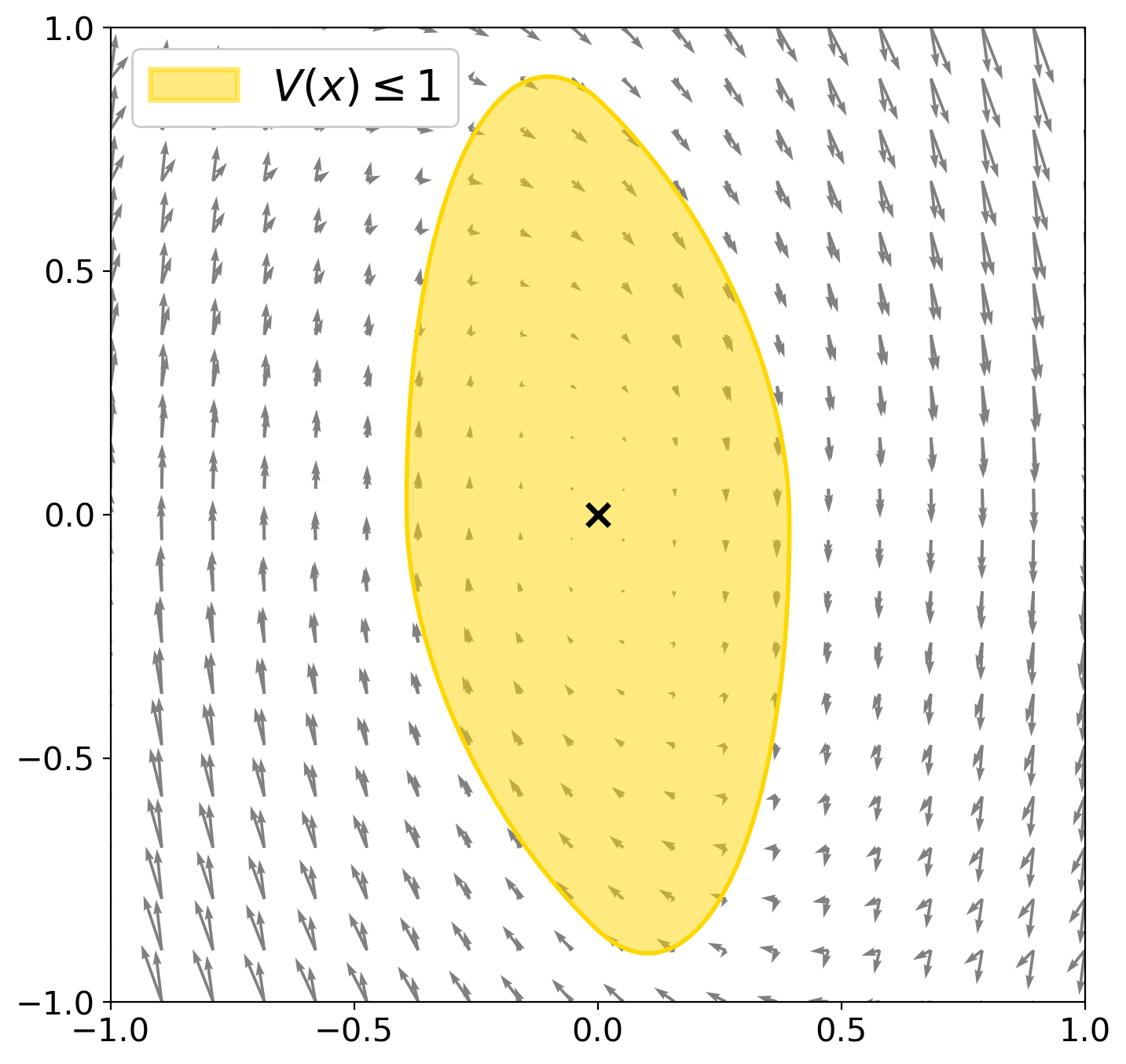}
\caption{Polyhedral Lyapunov function for the system described in
Subsection~\ref{ssec-benchmark-zelentsowsky} with $\alpha:6$ (gray arrows).
The function contains $300$ linear pieces.}
\label{fig-exa-zelentsowsky}
\end{figure}

Using Algorithm~\ref{algo-process}, we can compute a polyhedral Lyapunov function for the system with $\alpha:6$.
We use the parameters $\epsilon:50$, $\theta:1/64$ and $\delta:0.001$.
The computation takes about $5$ minutes, and the resulting function contains $300$ linear pieces (see Figure~\ref{fig-exa-zelentsowsky}).
We also show that the system with $\alpha:6.87$ does not admit a polyhedral Lyapunov function with these parameters.

\begin{rem}
Let us mention that the algorithms in the above papers focus on uncertain linear systems,
while our algorithm also tackles \emph{piecewise} uncertain linear systems.
Another difference, e.g., with~\citet{ambrosino2012aconvex}, is that our algorithm does not involve any hyper-parameters
that must be set by the user; it just requires the inputs $\epsilon$, $\theta$ and $\delta$,
that are used to verify whether the system admits a polyhedral Lyapunov function with specified eccentricity and robustness.
\end{rem}

%%%%%%%%%%%%%%%%%%%%%%%%%%%%%%%%%%%%%%%%%%%%%%%%%%%%%%%%%%%%%%%%%%%%%%%%%%%%%%%%%%%%%%%%%%%%%%%%%%%%
\subsection{Controlled mass--spring system}\label{ssec-mass-spring-example}

\begin{figure}
\centering
% \begin{subfigure}[c]{0.39\linewidth}
% 	\centering
% 	\includegraphics[width=\linewidth]{./figures/fig_exa_mass_spring_system}
% 	\caption{Mass-spring system}
% 	\label{fig-exa-mass-spring-system}
% \end{subfigure}
% \hfill
\begin{subfigure}[c]{\linewidth}
	\centering
	\includegraphics[width=0.9\linewidth]{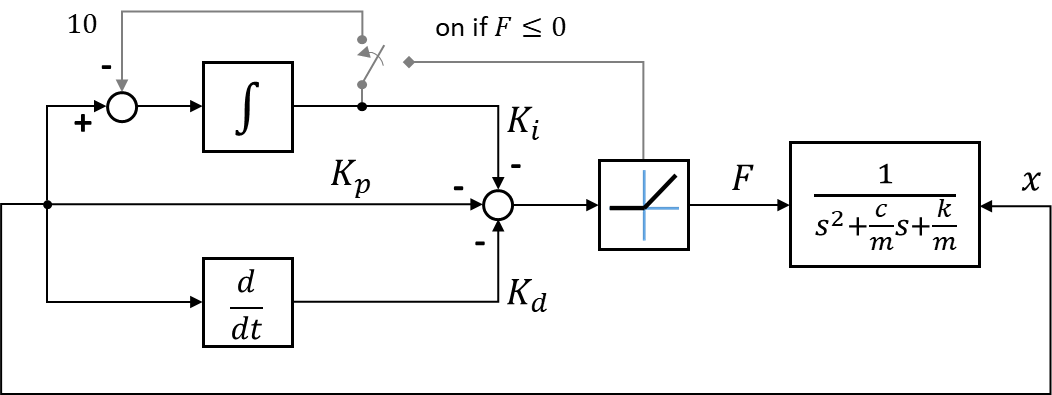}
	\caption{Block-diagram of the system system actuated by a ReLU-saturated PID controller.
	When the actuation force is saturated, an anti-windup mechanism (in gray)
	counterbalances the accumulation of the integrated error.}
	\label{fig-exa-mass-spring-diagram}
\end{subfigure}
\begin{subfigure}[c]{\linewidth}
	\centering
	\includegraphics[width=0.6\linewidth]{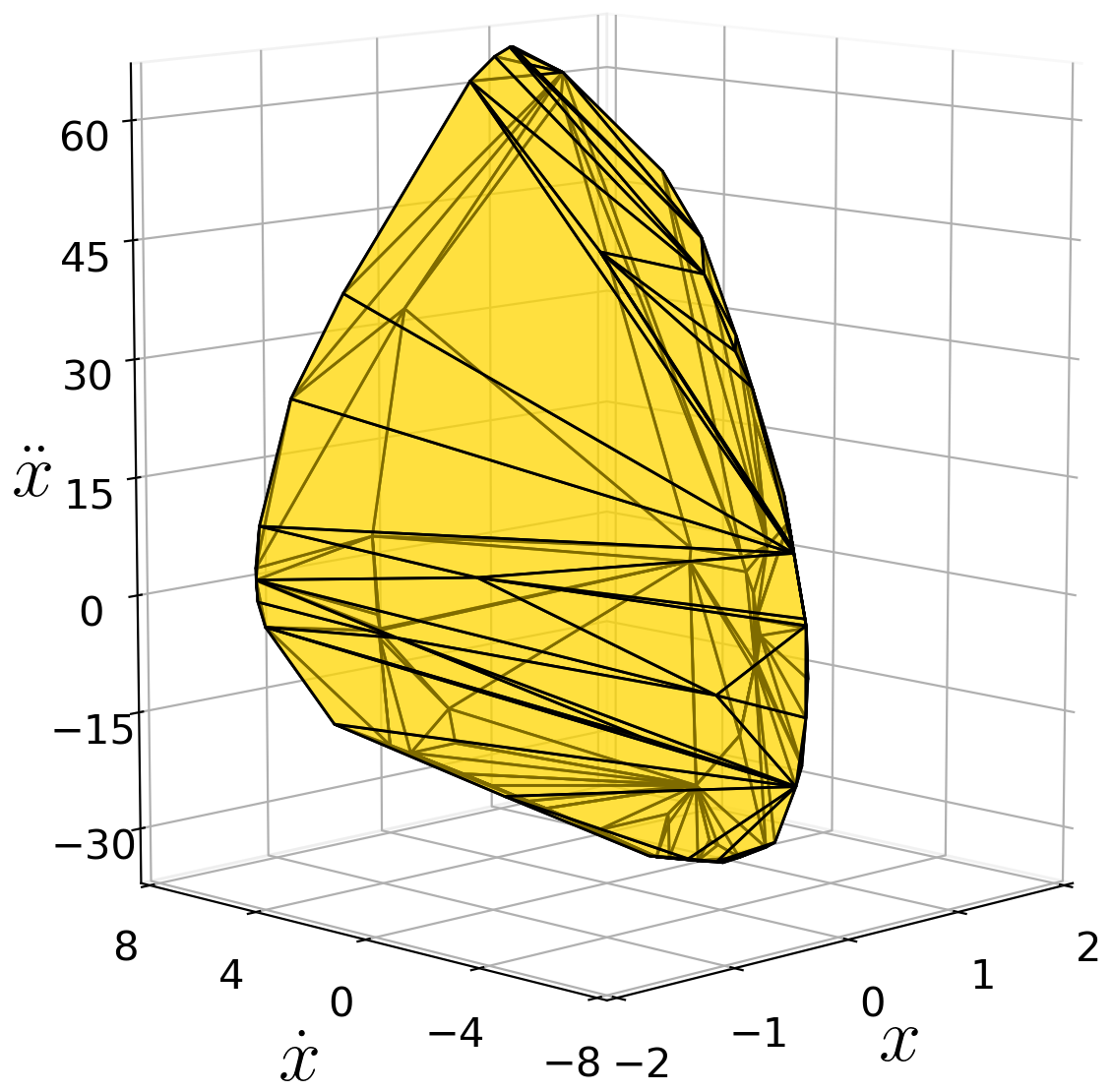}
	\caption{$1$-sublevel set of the Lyapunov function}
	\label{fig-exa-mass-spring-lyapunov}
\end{subfigure}
% \begin{subfigure}[c]{\linewidth}
% 	\centering
% 	\vskip3pt
% 	\includegraphics[width=\linewidth]{./figures/fig_exa_mass_spring_decrease}
% 	\vskip-7pt
% 	\caption{Value of the Lyapunov function along a trajectory}
% 	\label{fig-exa-mass-spring-decrease}
% \end{subfigure}
\caption{Example of Subsection~\ref{ssec-mass-spring-example}}
\label{fig-exa-mass-spring}
\end{figure}

We consider a mass--spring system whose dynamics is described by
$\ddot{x} = - \frac{k}{m}x + \frac1m F$,
where $m=0.1\,\mathrm{kg}$ and $k=2\,\Nrm/\mrm$.
We control this system using a PID controller defined by
$F(t) = - K_i y(t) - K_p x(t) - K_d \dot{x}(t)$,
where $y(t)=\int_0^t x(\tau)\,\diff\tau$, $K_i=44\,\Nrm/\mrm\cdot\srm$,
$K_p=24\,\Nrm/\mrm$ and $K_d=3.2\,\Nrm\cdot\srm/\mrm$.
The force that can be applied on the system can only be \emph{nonnegative}.
To counterbalance the accumulation of the error when the input is negative, we add an \emph{anti-windup} mechanism to the system.
The block-diagram of the resulting system is depicted in Figure~\ref{fig-exa-mass-spring-diagram}.
The dynamics of the system is described by the following piecewise linear system:
{\allowdisplaybreaks\begin{subequations}\label{eq-mass-spring-mode}
\begin{align}
&\begin{array}{l}
\text{if $K_d\dot{x} + K_px + K_iy \leq 0$:} \\
\hspace{1cm}\left\lbrace\begin{array}{l}
\dot{y} = x, \\
m\ddot{x} + K_d\dot{x} + (k+K_p)x + K_iy = 0, \\
\end{array}\right.
\end{array} \label{eq-mass-spring-mode1}\\
&\begin{array}{l}
\text{if $K_d\dot{x} + K_px + K_iy \geq 0$:} \\
\hspace{1cm}\left\lbrace\begin{array}{l}
\dot{y} = -y + x, \\
m\ddot{x} + kx = 0. \\
\end{array}\right.
\end{array} \label{eq-mass-spring-mode2}
\end{align}
\end{subequations}}%

\begin{rem}
Since~\eqref{eq-mass-spring-mode2} is not stable, System~\eqref{eq-mass-spring-mode}
does not admit a Lyapunov function symmetric around the origin
(including any polynomial Lyapunov function).
\end{rem}

Using Algorithm~\ref{algo-process}, we compute a polyhedral Lyapunov function for this system;
thereby showing that it is asymptotically stable.
We use the parameters $\epsilon:50$, $\theta:1/32$ and $\delta:0.001$.
The computation takes 15 seconds and outputs a polyhedral Lyapunov function
with $127$ linear pieces (see Figure~\ref{fig-exa-mass-spring-lyapunov}).
% We have also represented, in Figure~\ref{fig-exa-mass-spring-decrease}, the evolution of $V(x(t))$,
% where $x(t)$ is the trajectory of the system starting from $(y,x,\dot{x})=(0,1,0)$.
% We observe that $V(x(t))$ decreases with $t$, as predicted by the theory of Lyapunov.

%%%%%%%%%%%%%%%%%%%%%%%%%%%%%%%%%%%%%%%%%%%%%%%%%%%%%%%%%%%%%%%%%%%%%%%%%%%%%%%%%%%%%%%%%%%%%%%%%%%%
\subsection{Performance evaluation}\label{ssec-performance-example}

We want to evaluate the performance of the process,
in terms of computation time and complexity of the resulting Lyapunov function,
as a function of the dimension of the system and the stability margin of its matrices.
Therefore, for $d\in\Ne_{>0}$, we let $U\in\Re^{\dxd}$ be an orthogonal matrix,
and for a parameter $\gamma\in\Re$, we define the $\dxd$ matrix:
\[
\Pi_\gamma : U(\mathbf1\mathbf1^\top - (d+\gamma)I)U^\top,\quad \mathbf1=[1,\ldots,1]^\top\!\in\Re^d.
\]
Consider System $\calF$ with $Q:\{1,2\}$, $H_1:\Re_{\geq0}\times\Re^{d-1}$
and $H_2:\Re_{\leq0}\times\Re^{d-1}$, $A_1:\Pi_1$ and $A_2:\Pi_\gamma$.
For $\gamma<0$, the system is unstable (thus, does not admit any Lyapunov function);
see Lemma~\ref{lem-exa-performance-unstable} in Appendix~\ref{app-sec-results-ssec-performance-example}.
For $\gamma>0$, the system is asymptotically stable and admits a $2d$-piece polyhedral Lyapunov function;
see Lemma~\ref{lem-exa-performance-stable} in Appendix~\ref{app-sec-results-ssec-performance-example}.
% It follows that, for $\gamma_2\leq\gamma$, any $(0,\gamma_2)$-perturbation
% (Definition~\ref{def-perturbed-system}) of System $\calF$ is stable;
% reversely if $\gamma_2>\gamma$, then there exists a $(0,\gamma_2)$-perturbation of System $\calF$ that is unstable.

For different values of $d\in\Ne_{>0}$ and $\gamma\in\Re_{>0}$,
we consider the parameters $\epsilon:10$, $\theta:1/8$ and $\delta\in\{\frac\gamma{50},\gamma\}$.
We use Algorithm~\ref{algo-process} to compute a polyhedral Lyapunov function for System $\calF$,
or conclude that no polyhedral Lyapunov function with the given parameters exists.
For each case, we measure the computation time and the number of linear pieces of the obtained Lyapunov function,
or the number of iterations before the algorithm outputs \textsc{fail}.
The results are gathered in Table~\ref{tab-performance}.
We observe that the computation time and the number of pieces/iterations increase
when $d$ increases and/or $\gamma$ decreases.
This is in accordance with the complexity analysis given in Corollary~\ref{cor-termination-process},
which state that the maximal number of iterations depends on the $\rbar$-packing number in dimension $d$,
wherein $\rbar$ depends, among others, on $\gamma$.

\begin{table}
\centering
\renewcommand{\arraystretch}{1.1}
\begin{tabular}{@{}llc@{}ll@{}c@{}ll@{}}
\toprule
&&\phantom{abcd}& \multicolumn{2}{@{}l}{$\delta = \frac\gamma{50}$} &\phantom{abc}& \multicolumn{2}{@{}l}{$\delta = \gamma$} \\
\cmidrule{4-5}\cmidrule{7-8}
& $\gamma$ && $T$ [sec] & pieces && $T$ [sec] & iterations \\
\midrule
$d = 4$~~
& 1 && 0.37 & 12 && 0.01 & 2\\
& 0.1 && 0.84 & 21 && 0.23 & 13\\
& 0.01 && 4.48 & 37 && 1.17 & 25\\
\midrule
$d = 5$
& 1 && 0.47 & 14 && 0.01 & 2\\
& 0.1 && 2.82 & 30 && 0.33 & 14\\
& 0.01 && 11.54 & 53 && 2.47 & 31\\
\midrule
$d = 6$
& 1 && 1.31 & 20 && 0.01 & 2\\
& 0.1 && 11.67 & 46 && 1.07 & 22\\
& 0.01 && 12.32 & 52 && 4.58 & 40\\
\midrule
$d = 7$
& 1 && 2.22 & 23 && 0.01 & 2\\
& 0.1 && 26.42 & 61 && 2.16 & 25\\
& 0.01 && 34.09 & 68 && 11.94 & 41\\
\midrule
$d = 8$
& 1 && 20.20 & 50 && 0.01 & 2\\
& 0.1 && 47.27 & 87 && 5.65 & 35\\
& 0.01 && 160.18 & 133 && 7.95 & 46\\
\midrule
$d = 9$
& 1 && 13.85 & 45 && 0.01 & 2\\
& 0.1 && 92.77 & 98 && 10.22 & 42\\
& 0.01 && 206.85 & 149 && 24.61 & 54\\
\bottomrule
\end{tabular}
\vskip5pt
\caption{Evaluation of the performance of Algorithm~\ref{algo-process} on
System $\calF$ described in Subsection~\ref{ssec-performance-example}.
``$T$'' refers to the total computation time of the algorithm.
``Pieces'' refers to the number of linear pieces of the computed polyhedral Lyapunov function,
and ``iterations'' refers to the number of iterations needed by the algorithm to
conclude that no polyhedral Lyapunov function exists.}
\label{tab-performance}
\end{table}

%%%%%%%%%%%%%%%%%%%%%%%%%%%%%%%%%%%%%%%%%%%%%%%%%%%%%%%%%%%%%%%%%%%%%%%%%%%%%%%%%%%%%%%%%%%%%%%%%%%%
\section{Conclusions}

We presented an algorithmic framework to compute polyhedral
Lyapunov functions for continuous-time piecewise linear systems.
Compared to previous approaches in the literature, a key asset of our approach
is that we do not put a priori bound on the number of facets of the Lyapunov function
and we provide formal guarantees on the computational complexity of the algorithm.
While these theoretical guarantees can be large when the dimension of the system increases,
we demonstrate their practical applicability on numerical examples of dimensions from $2$ to $9$.

For further work, we will consider the problem of controller synthesis for linear hybrid systems,
by extending the approach for the learning of polyhedral \emph{control} Lyapunov functions (CLFs)
from counterexamples.
This aim to address a lot of interesting control problems where the dynamics of a complex system
(e.g., robot dynamics, power electronics) can be abstracted by linear hybrid systems (such as Neural Networks).

 %% ADDED BY GUILLAUME

\begin{ack}                               % Place acknowledgements
To be added in the final version. 
\end{ack}

\bibliography{myrefs}

\begin{thebibliography}{30}
\providecommand{\natexlab}[1]{#1}
\providecommand{\url}[1]{\texttt{#1}}
\expandafter\ifx\csname urlstyle\endcsname\relax
  \providecommand{\doi}[1]{doi: #1}\else
  \providecommand{\doi}{doi: \begingroup \urlstyle{rm}\Url}\fi

\bibitem[Abate et~al.(2021)Abate, Ahmed, Giacobbe, and
  Peruffo]{abate2021formal}
Alessandro Abate, Daniele Ahmed, Mirco Giacobbe, and Andrea Peruffo.
\newblock Formal synthesis of {Lyapunov} neural networks.
\newblock \emph{IEEE Control Systems Letters}, 5\penalty0 (3):\penalty0
  773--778, 2021.
\newblock \doi{10.1109/LCSYS.2020.3005328}.

\bibitem[Ahmadi and Jungers(2016)]{ahmadi2016lower}
Amir~Ali Ahmadi and Rapha{\"e}l~M Jungers.
\newblock Lower bounds on complexity of {Lyapunov} functions for switched
  linear systems.
\newblock \emph{Nonlinear Analysis: Hybrid Systems}, 21:\penalty0 118--129,
  2016.
\newblock \doi{10.1016/j.nahs.2016.01.003}.

\bibitem[Ambrosino et~al.(2012)Ambrosino, Ariola, and
  Amato]{ambrosino2012aconvex}
Roberto Ambrosino, Marco Ariola, and Francesco Amato.
\newblock A convex condition for robust stability analysis via polyhedral
  {Lyapunov} functions.
\newblock \emph{SIAM Journal on Control and Optimization}, 50\penalty0
  (1):\penalty0 490--506, 2012.
\newblock \doi{10.1137/110796479}.

\bibitem[Berger and Sankaranarayanan(2022)]{berger2022learning}
Guillaume~O Berger and Sriram Sankaranarayanan.
\newblock Learning fixed-complexity polyhedral {Lyapunov} functions from
  counterexamples, 2022.
\newblock \arxiv{2204.06693} (v1: 14\textsuperscript{th} April 2022).

\bibitem[Blanchini and Miani(1996)]{blanchini1996onthetransient}
Franco Blanchini and Stefano Miani.
\newblock On the transient estimate for linear systems with time-varying
  uncertain parameters.
\newblock \emph{IEEE Transactions on Circuits and Systems I: Fundamental Theory
  and Applications}, 43\penalty0 (7):\penalty0 592--596, 1996.
\newblock \doi{10.1109/81.508180}.

\bibitem[Blanchini and Miani(2015)]{blanchini2015settheoretic}
Franco Blanchini and Stefano Miani.
\newblock \emph{Set-theoretic methods in control}.
\newblock Birkh{\"a}user, Cham, 2\textsuperscript{nd} edition, 2015.
\newblock \doi{10.1007/978-3-319-17933-9}.

\bibitem[Boyd and Vandenberghe(2004)]{boyd2004convex}
Stephen Boyd and Lieven Vandenberghe.
\newblock \emph{Convex optimization}.
\newblock Cambridge University Press, Cambridge, UK, 2004.
\newblock \doi{10.1017/CBO9780511804441}.

\bibitem[Chang et~al.(2019)Chang, Roohi, and Gao]{chang2019neural}
Ya-Chien Chang, Nima Roohi, and Sicun Gao.
\newblock Neural {Lyapunov} control.
\newblock In \emph{NIPS'19: Proceedings of the 33rd International Conference on
  Neural Information Processing Systems}, pages 3245--3254. ACM, 2019.
\newblock \doi{10.5555/3454287.3454579}.

\bibitem[Chesi et~al.(2009)Chesi, Garulli, Tesi, and
  Vicino]{chesi2009homogeneous}
Graziano Chesi, Andrea Garulli, Alberto Tesi, and Antonio Vicino.
\newblock \emph{Homogeneous polynomial forms for robustness analysis of
  uncertain systems}, volume 390 of \emph{Lecture Notes in Control and
  Information Sciences}.
\newblock Springer, London, 2009.
\newblock \doi{10.1007/978-1-84882-781-3}.

\bibitem[Dai et~al.(2021)Dai, Landry, Yang, Pavone, and
  Tedrake]{dai2021lyapunovstable}
Hongkai Dai, Benoit Landry, Lujie Yang, Marco Pavone, and Russ Tedrake.
\newblock {Lyapunov-stable} neural-network control.
\newblock In \emph{Proceedings of Robotics: Science and Systems}, Virtual,
  2021.
\newblock \doi{10.15607/RSS.2021.XVII.063}.

\bibitem[Goebel et~al.(2012)Goebel, Sanfelice, and Teel]{goebel2012hybrid}
Rafal Goebel, Ricardo~G Sanfelice, and Andrew~R Teel.
\newblock \emph{Hybrid dynamical systems: modeling stability, and robustness}.
\newblock Princeton University Press, Princeton, NJ, 2012.

\bibitem[Guglielmi et~al.(2017)Guglielmi, Laglia, and
  Protasov]{guglielmi2017polytope}
Nicola Guglielmi, Linda Laglia, and Vladimir Protasov.
\newblock Polytope {Lyapunov} functions for stable and for stabilizable {LSS}.
\newblock \emph{Foundations of Computational Mathematics}, 17\penalty0
  (2):\penalty0 567--623, 2017.
\newblock \doi{10.1007/s10208-015-9301-9}.

\bibitem[Hassibi and Boyd(1998)]{hassibi1998quadratic}
Arash Hassibi and Stephen Boyd.
\newblock Quadratic stabilization and control of piecewise-linear systems.
\newblock In \emph{Proceedings of the 1998 American Control Conference. ACC
  (IEEE Cat. No. 98CH36207)}, pages 3659--3664. IEEE, 1998.
\newblock \doi{10.1109/ACC.1998.703296}.

\bibitem[Jungers(2009)]{jungers2009thejoint}
Rapha{\"e}l~M Jungers.
\newblock \emph{The joint spectral radius: theory and applications}.
\newblock Springer, Berlin, 2009.
\newblock \doi{10.1007/978-3-540-95980-9}.

\bibitem[Kapinski et~al.(2014)Kapinski, Deshmukh, Sankaranarayanan, and
  Ar{\'e}chiga]{kapinski2014simulationguided}
James Kapinski, Jyotirmoy~V Deshmukh, Sriram Sankaranarayanan, and Niko
  Ar{\'e}chiga.
\newblock Simulation-guided {Lyapunov} analysis for hybrid dynamical systems.
\newblock In \emph{Proceedings of the 17th international conference on Hybrid
  systems: computation and control}, pages 133--142. ACM, 2014.
\newblock \doi{10.1145/2562059.2562139}.

\bibitem[Khalil(2002)]{khalil2002nonlinear}
Hassan~K Khalil.
\newblock \emph{Nonlinear systems}.
\newblock Prentice-Hall, Upper Saddle River, NJ, 3\textsuperscript{rd} edition,
  2002.

\bibitem[Kousoulidis and Forni(2021)]{kousoulidis2021polyhedral}
Dimitris Kousoulidis and Fulvio Forni.
\newblock Polyhedral {Lyapunov} functions with fixed complexity.
\newblock In \emph{2021 60th IEEE Conference on Decision and Control (CDC)},
  pages 3293--3298. IEEE, 2021.
\newblock \doi{10.1109/CDC45484.2021.9683697}.

\bibitem[Lasserre(2001)]{lasserre2001global}
Jean~B Lasserre.
\newblock Global optimization with polynomials and the problem of moments.
\newblock \emph{SIAM Journal on Optimization}, 11\penalty0 (3):\penalty0
  796--817, 2001.
\newblock \doi{10.1137/S1052623400366802}.

\bibitem[Lazar and Doban(2011)]{lazar2011oninfinity}
Mircea Lazar and Alina~I Doban.
\newblock On infinity norms as {Lyapunov} functions for continuous-time
  dynamical systems.
\newblock In \emph{2011 50th IEEE Conference on Decision and Control and
  European Control Conference}, pages 7567--7572. IEEE, 2011.
\newblock \doi{10.1109/CDC.2011.6161163}.

\bibitem[Liberzon(2003)]{liberzon2003switching}
Daniel Liberzon.
\newblock \emph{Switching in systems and control}.
\newblock Birkh{\"a}user, Boston, MA, 2003.
\newblock \doi{10.1007/978-1-4612-0017-8}.

\bibitem[Miani and Savorgnan(2005)]{miani2005maxisg}
Stefano Miani and Carlo Savorgnan.
\newblock {MAXIS-G}: a software package for computing polyhedral invariant sets
  for constrained {LPV} systems.
\newblock In \emph{Proceedings of the 44th IEEE Conference on Decision and
  Control}, pages 7609--7614. IEEE, 2005.
\newblock \doi{10.1109/CDC.2005.1583390}.

\bibitem[Nesterov and Nemirovskii(1994)]{nesterov1994interiorpoint}
Yurii Nesterov and Arkadii Nemirovskii.
\newblock \emph{Interior-point polynomial algorithms in convex programming}.
\newblock SIAM, Philadelphia, PA, 1994.
\newblock \doi{10.1137/1.9781611970791}.

\bibitem[Pettit and Wellstead(1995)]{pettit1995analyzing}
Njal~BOL Pettit and Peter~E Wellstead.
\newblock Analyzing piecewise linear dynamical systems.
\newblock \emph{IEEE Control Systems Magazine}, 15\penalty0 (5):\penalty0
  43--50, 1995.
\newblock \doi{10.1109/37.466263}.

\bibitem[Pola{\'n}ski(2000)]{polanski2000onabsolute}
Andrzej Pola{\'n}ski.
\newblock On absolute stability analysis by polyhedral {Lyapunov} functions.
\newblock \emph{Automatica}, 36\penalty0 (4):\penalty0 573--578, 2000.
\newblock \doi{10.1016/S0005-1098(99)00180-6}.

\bibitem[Prabhakar and Viswanathan(2013)]{prabhakar2013onthedecidability}
Pavithra Prabhakar and Mahesh Viswanathan.
\newblock On the decidability of stability of hybrid systems.
\newblock In \emph{HSCC'13: Proceedings of the 16th international conference on
  Hybrid systems: computation and control}, pages 53--62. ACM, 2013.
\newblock \doi{10.1145/2461328.2461339}.

\bibitem[Ravanbakhsh and Sankaranarayanan(2019)]{ravanbakhsh2019learning}
Hadi Ravanbakhsh and Sriram Sankaranarayanan.
\newblock Learning control {Lyapunov} functions from counterexamples and
  demonstrations.
\newblock \emph{Autonomous Robots}, 43\penalty0 (2):\penalty0 275--307, 2019.
\newblock \doi{10.1007/s10514-018-9791-9}.

\bibitem[Sun and Ge(2011)]{sun2011stability}
Zhendong Sun and Shuzhi~Sam Ge.
\newblock \emph{Stability theory of switched dynamical systems}.
\newblock Springer, London, 2011.
\newblock \doi{10.1007/978-0-85729-256-8}.

\bibitem[Topcu et~al.(2008)Topcu, Packard, and Seiler]{topcu2008local}
Ufuk Topcu, Andrew Packard, and Peter Seiler.
\newblock Local stability analysis using simulations and sum-of-squares
  programming.
\newblock \emph{Automatica}, 44\penalty0 (10):\penalty0 2669--2675, 2008.
\newblock \doi{10.1016/j.automatica.2008.03.010}.

\bibitem[Xie et~al.(1997)Xie, Shishkin, and Fu]{xie1997piecewise}
Lin Xie, Serge Shishkin, and Minyue Fu.
\newblock Piecewise {Lyapunov} functions for robust stability of linear
  time-varying systems.
\newblock \emph{Systems \& Control Letters}, 31\penalty0 (3):\penalty0
  165--171, 1997.
\newblock \doi{10.1016/S0167-6911(97)00027-3}.

\bibitem[Zelentsovsky(1994)]{zelentsovsky1994nonquadratic}
AL~Zelentsovsky.
\newblock Nonquadratic {Lyapunov} functions for robust stability analysis of
  linear uncertain systems.
\newblock \emph{IEEE Transactions on Automatic Control}, 39\penalty0
  (1):\penalty0 135--138, 1994.
\newblock \doi{10.1109/9.273350}.

\end{thebibliography}

\appendix   % Each appendix must have a short title.
            % Sections and subsections are supported
            % in the appendices.

%%%%%%%%%%%%%%%%%%%%%%%%%%%%%%%%%%%%%%%%%%%%%%%%%%%%%%%%%%%%%%%%%%%%%%%%%%%%%%%%%%%%%%%%%%%%%%%%%%%%
\section{Proof of Theorem~\ref{thm-equivalent-piecewise-linear-lyapunov}}\label{app-sec-proof-thm-equivalent-piecewise-linear-lyapunov}

{\itshape Proof of ``if'' direction.}
It is straightforward to see that (D1)--(D2) implies
(C1)--(C2) in Proposition~\ref{pro-piecewise-linear-lyapunov}.
Therefore, by Proposition~\ref{pro-piecewise-linear-lyapunov},
(D1)--(D2) implies that $V$ is a Lyapunov function System $\calF$.

{\itshape Proof of ``only if'' direction.}
Assume that $V$ is a Lyapunov function for System $\calF$.
Then, by (C1) in Proposition~\ref{pro-piecewise-linear-lyapunov} and $V$ being continuous,
there is $\epsilon\geq1$ s.t.~for all $x\in\Sbb$, $V(x)\geq\frac\Vmax\epsilon$.
Now, since $V$ is positively homogeneous of degree $1$, it follows that for all $x\in\Re^d$,
$V(x)\geq\frac\Vmax\epsilon\lVert x\rVert$, so that (D1) holds.

It remains to find $\theta$ and $\delta$ s.t.~(D2) holds.
Therefore, fix $x\in\Sbb$ and $c\in\calV$.
First, assume that $c\in\calV(x)$.
Then, by (C2), there is $\delta_{x,c}>0$
s.t.~for all $v\in\calV$, $c^\top v<-\delta_{x,c}\Vmax$.
Now, assumme that $c\notin\calV(x)$, i.e., $c^\top x<V(x)$.
Thus, there is $\delta_{x,c}>0$ and $\theta_{x,c}>0$
s.t.~for all $v\in\calV$, $c^\top(x+\theta_{x,c}v)<V(x)-\delta_{x,c}\Vmax$.
Since $\calF$ and $V$ are continuous, there is a neighborhood $\calN_x\subseteq\Sbb$ of $x$
s.t.~for all $x'\in\calN_x$ and all $v\in\calF(x')$,
$c^\top(x'+\theta_{x,c}v)<V(x')-\delta_{x,c}\Vmax$.

Since $x$ was arbitrary, the above holds for every $x\in\Sbb$.
By compactness of $\Sbb$, there exists a finite set of points $x_1,\ldots,x_n\in\Sbb$
s.t.~$\bigcup_{i=1}^n\calN_{x_i}$ covers $\Sbb$.
We can then define $\theta:\min_{i:1,\ldots,n}\,\theta_{x_i,c}>0$
and $\delta:\min_{i:1,\ldots,n}\,\theta_{x_i,c}>0$.
(D2) is then satisfied for $c$ and for all $x\in\Sbb$.
By the positive homogeneity of $V$, it then follows that (D2) is satisfied for $c$ and for all $x\in\Re^d$.
Now, since there is a finite set of values for $c$,
we can find $\theta$ and $\delta$ so that (D2) holds, concluding the proof.\qed

\section{Proof of Theorem~\ref{thm-robustness-necessary}}\label{app-sec-proof-thm-robustness-necessary}

First, note that by (D1) in Definition~\ref{def-robust-lyapunov-conditions},
it holds that for all $x\in\Re^d\setminus\{0\}$ and $c\in\calV(x)$, $\lVert c\rVert_*\geq\frac\Vmax\epsilon$.
W.l.o.g.~assume that $\Vmax=1$.
Denote $\eta:\frac{\gamma\amax}\epsilon$.
We will need the following result.

\begin{lem}\label{lem-robust-c-gamma}
For all $x\in\Re^d$, $v\in\calF(x)$ and $c\in\calV(x)$, $c^\top v \leq -\eta\lVert x\rVert$.
\end{lem}

\begin{pf}
Let $x\in\Re^d$, $v\in\calF(x)$ and $c\in\calV(x)$
Since $V$ is a Lyapunov function for any $\gamma$-perturbation of System $\calF$,
it holds that for all $u\in\Re^d$ s.t.~$\lVert u\rVert\leq\gamma\amax\lVert x\rVert$, $c^\top(v+u)\leq0$.
This implies that $c^\top v\leq-\gamma\amax\lVert c\rVert_*\lVert x\rVert\leq -\frac{\gamma\amax}\epsilon\lVert x\rVert = -\eta\lVert x\rVert$.\qed
\end{pf}

We proceed with the proof of the theorem.
Therefore, fix $x\in\Re^d$, $q\in Q$ and $c\in\calV$.
Let $v:A_qx$.
We will show that $c^\top v\leq\frac1\theta(V(x)-c^\top x)-\frac\eta2\lVert x\rVert$, where $\theta:\frac\eta{2\amax^2+2\eta\amax}$.
First, if $c\in\calV(x)$, then we are done by Lemma~\ref{lem-robust-c-gamma}.
Thus, assume that $V(x)>c^\top x$.

Let $s:\frac2{\eta\lVert x\rVert}(V(x)-c^\top x)>0$ and $y:x+sv$.
Let $c_y\in\calV(y)$.
We will need the following result.

\begin{lem}\label{lem-cv-cyv}
$c^\top v \leq c_y^\top v + \frac\eta2\lVert x\rVert$.
\end{lem}

\begin{pf}
$sc_y^\top v \geq c_y^\top x - V(x) + sc_y^\top v = c_y^\top y - V(x)
\geq c^\top y - V(x) = c^\top x + sc^\top v - V(x) = s(c^\top v - \frac\eta2\lVert x\rVert)$.\qed
\end{pf}

Note that $\lVert v\rVert\leq\amax\lVert x\rVert$, so that $\lVert y-x\rVert\leq s\amax\lVert x\rVert$.
Let $w:A_qy$.
By Lemma~\ref{lem-robust-c-gamma}, it holds that $c_yw\leq-\eta\lVert y\rVert$.
Thus, $c_y^\top w\leq-\eta(1-s\amax)\lVert x\rVert$.
Also, $\lVert w-v\rVert=\lVert A_q(y-x)\rVert\leq s\amax^2\lVert x\rVert$.
Hence, by Lemma~\ref{lem-cv-cyv}, $c^\top v\leq s(\amax^2+\eta\amax)\lVert x\rVert-\frac\eta2\lVert x\rVert$.
From the definition of $s$, this gives $c^\top v\leq\frac1\theta(V(x)-c^\top x)-\frac\eta2\lVert x\rVert$.
Since $x$, $v$ and $c$ were arbitrary, this concludes the proof.\qed

%%%%%%%%%%%%%%%%%%%%%%%%%%%%%%%%%%%%%%%%%%%%%%%%%%%%%%%%%%%%%%%%%%%%%%%%%%%%%%%%%%%%%%%%%%%%%%%%%%%%
\section{Results in Subsection~\ref{ssec-performance-example}}\label{app-sec-results-ssec-performance-example}

\begin{lem}\label{lem-exa-performance-unstable}
For $\gamma<0$, the system described in Subsection~\ref{ssec-performance-example} is unstable.
\end{lem}

\begin{pf}
Let $x:\pm U\mathbf1$, where the sign  ``$\pm$'' is chosen is such a way that $x\in H_2$.
It holds that $x$ is an eigenvector of $\Pi_\gamma$ with eigenvalue $-\gamma>0$,
so that the trajectory starting from $x$ is divergent.\qed
\end{pf}

\begin{lem}\label{lem-exa-performance-stable}
For $\gamma>0$, the system described in Subsection~\ref{ssec-performance-example}
admits a $2d$-piece polyhedral Lyapunov function, e.g.,
the polyhedral function $x\mapsto\lVert U^\top x\rVert_\infty$.
\end{lem}

\begin{pf}
W.l.o.g.~(using a change of coordinates if necessary) assume that $U=I$.
Let $V:\lVert\cdot\rVert_\infty$.
Since $\lVert\cdot\rVert_\infty$ is a norm,
it is clear that $V$ satisfies (C1) in Proposition~\ref{pro-piecewise-linear-lyapunov}.
To show that $V$ satisfies (C2) in Proposition~\ref{pro-piecewise-linear-lyapunov},
let $x\in\Re^d$.
Let $e\in\{0,1\}^d$ and $\sigma\in\{-1,1\}$ be s.t.~$e^\top\mathbf1=1$ and $V(x)=\sigma e^\top x$.
It holds that $\sigma e^\top\Pi_\gamma x=\sigma\mathbf1^\top x-\sigma(d+\gammabar)e^\top x\leq-\sigma\gammabar e^\top x=-\gammabar V(x)$,
where $\gammabar\in\{1,\gamma\}$ depending on whether $x\in H_1$ or $H_2$.
Hence, $V$ satisfies (C2) in Proposition~\ref{pro-piecewise-linear-lyapunov}, concluding the proof.\qed
\end{pf} %% ADDED BY GUILLAUME

\end{document}